\documentclass[11pt,a4paper]{article}

\usepackage{mathtools}
\usepackage{authblk} 
\usepackage{algpseudocode,algorithmicx,algorithm}

\usepackage{mathrsfs}
\usepackage{latexsym,bm,bbm}
\usepackage{amsmath,amsfonts,amssymb,amsthm,rotating}
\usepackage{extarrows}

\usepackage{graphicx,subfigure,epstopdf,float}
\usepackage{enumerate,cases,multirow}
\usepackage{makecell}
\usepackage{caption}

\usepackage{longtable,colortbl,arydshln,threeparttable}
\definecolor{mygray}{gray}{.9}

\usepackage{indentfirst}
\usepackage[top=25mm,bottom=20mm,left=25mm,right=20mm]{geometry}
\usepackage{cite}

\usepackage{listings}

\usepackage{makeidx}        
\usepackage{booktabs}
\usepackage{xcolor}
\usepackage[bookmarksnumbered,colorlinks=true,citecolor=red,linkcolor=red,hyperindex,linktocpage=true]{hyperref}

\newcommand{\ket}[1]{| #1 \rangle} 
\newcommand{\bra}[1]{\langle #1 |} 

\newcommand{\ketbra}[2]{ | #1 \rangle \langle #2 |}
\newcommand{\bb}{\boldsymbol}

\def \d {\mathrm{d}}
\def \e {\mathrm{e}}
\def \i {\mathrm{i}}

\newcounter{parentalgorithm}

\makeatother

\newtheorem{theorem}{Theorem}[section]
\newtheorem{lemma}{Lemma}[section]
\newtheorem{proposition}{Proposition}[section]

\newtheorem{definition}{Definition}[section]

\theoremstyle{remark}
\newtheorem{remark}{\bf Remark}[section]

\numberwithin{equation}{section}

\begin{document}

\title{Schr\"odingerization based quantum algorithms for the time-fractional heat equation}

\author[1]{Shi Jin\thanks{shijin-m@sjtu.edu.cn}}
\author[1, 2]{Nana Liu\thanks{nana.liu@quantumlah.org}}
\author[3]{Yue Yu\thanks{terenceyuyue@xtu.edu.cn}}
\affil[1]{School of Mathematical Sciences, Institute of Natural Sciences, MOE-LSC, Shanghai Jiao Tong University, Shanghai, 200240, China}
\affil[2]{Global College, Shanghai Jiao Tong University, Shanghai 200240, China}
\affil[3]{School of Mathematics and Computational Science, Hunan Key Laboratory for Computation and Simulation in Science and Engineering, Key Laboratory of Intelligent Computing and Information Processing of Ministry of Education, National Center for Applied Mathematics in Hunan, Xiangtan University, Xiangtan, Hunan 411105, China}

\maketitle

\begin{abstract}
  We develop a quantum algorithm for solving high-dimensional time-fractional heat equations. By applying the dimension extension technique from \cite{Eritz23timeFractional}, the $d+1$-dimensional time-fractional equation is reformulated as a local partial differential equation in $d+2$ dimensions. Through discretization along both the extended and spatial domains,  a stable system of ordinary differential equations is obtained  by a simple change of variables. We propose a quantum algorithm for the resulting semi-discrete problem using the Schr\"odingerization approach from \cite{JLY22SchrShort, JLY22SchrLong, analogPDE}. The Schr\"odingerization technique transforms general linear partial and ordinary differential equations into Schr\"odinger-type systems~--~with unitary evolution, making them suitable for quantum simulation. This is accomplished via the warped phase transformation, which maps the equation into a higher-dimensional space. We provide detailed implementations of this method and conduct a comprehensive complexity analysis, demonstrating up to exponential advantage~--~with respect to the inverse of the mesh size in high dimensions~--~compared to its classical counterparts. Specifically, to compute the solution to time $T$,  while the classical method requires at least $\mathcal{O}(N_t d h^{-(d+0.5)})$ matrix-vector multiplications, where $N_t $ is the number of time steps (which is, for example, $ \mathcal{O}(Tdh^{-2})$ for the forward Euler method), our quantum algorithms requires $\widetilde{\mathcal{O}}(T^2d^4 h^{-8})$ queries to the block-encoding input models, with the quantum complexity being independent of the dimension $d$ in terms of the inverse mesh size $h^{-1}$.  Numerical experiments are performed to validate our formulation.
\end{abstract}



\section{Introduction}

Fractional differential equations are important  in models of complex physical systems, as they offer more realistic characterizations of certain phenomena compared to traditional differential operators. In particular, time-fractional differential equations are used in applications such as modeling tumor growth, anomalous diffusion, and certain viscoelastic models \cite{Khristenko23fractional,Eritz23timeFractional,Podlubny1995fracheat,Mainardi2010fracelasticity}. Given the challenges in obtaining analytical solutions for these problems, the development of numerical methods has become essential. For a brief review of classical methods for fractional ordinary differential equations (ODEs), we refer the reader to \cite{Khristenko23fractional} and the references therein.

Quantum computing offers a promising approach to overcome the computational challenges associated with high dimensionality, as it requires only a logarithmic number of qubits relative to the matrix dimension to store both the matrix and the solution vector. This advantage has led to increasing interest in developing quantum algorithms for solving PDEs \cite{Cao2013Poisson,Berry2014Highorder,qFEM-2016,Costa2019Wave,
Engel2019qVlasov,Childs2020qSpectral,Linden2020heat,Childs2021high,JinLiu2022nonlinear,GJL2022QuantumUQ,JLY2022multiscale}. Among the notable quantum strategies for solving PDEs is the {\it Schr\"odingerization} method introduced in \cite{JLY22SchrShort, JLY22SchrLong, analogPDE}. This simple and generic framework enables quantum simulation for {\it all} linear PDEs and ODEs. The core idea is to apply a warped phase transformation that maps the equations to one higher dimension, which, in the Fourier space, transforms them into a system of Schr\"odinger-type equations~--~with {\it unitary} evolutions!

The approach has been expanded to address a wide array of problems, including open quantum systems within bounded domains with non-unitary artificial boundary conditions \cite{JLLY23ABC}, problems entailing physical boundary or interface conditions \cite{JLLY2024boundary,JLY24Circuits}, Maxwell's equations \cite{JLC23Maxwell,MJL2024time}, the Fokker-Planck equation \cite{JLY24FokkerPlanck}, multiscale PDEs \cite{HJZ2024multiscale}, stochastic differential equations \cite{JLW2024Stochastic}, ill-posed problems such as the backward heat equation \cite{JLM24SchrBackward}, linear dynamical systems with inhomogeneous terms \cite{JLM24SchrInhom}, non-autonomous ordinary and partial differential equations \cite{CJL23TimeSchr,CJL2024TimeSchr}, iterative linear algebra solvers \cite{JinLiu-LA}, highly-oscillatory transport equations \cite{GJ2025transport}, etc. This approach is also natural for continuous variables quantum computing \cite{analogPDE} thus provides a framework~--~so far the only possible one~--~for the analog quantum simulation of PDEs and ODEs that completely retains the continuous nature of the system  without numerical discretizations. A complementary approach based on integral representation of the solution followed by Linear Combination of Hamiltonian Simulations was given in \cite{ALL2023LCH}, whose relation to Schr\"odingerization is studied in \cite{Li25}.

 As demonstrated in the Schr\"odingerization approach, dimension lifting has proven to be a versatile method applicable in the design of quantum algorithms for linear differential equations. This perspective extends to various other contexts, such as transforming nonlinear PDEs into linear ones \cite{JinLiu2022nonlinear}, converting non-autonomous Hamiltonian systems with time-dependent Hamiltonians into autonomous PDEs with time-independent coefficients \cite{CJL23TimeSchr,CJL2024TimeSchr,JLY24Circuits,MJL2024time}, and developing new transformations that map PDEs with uncertainty to deterministic PDEs \cite{GJL2022QuantumUQ}. While dimension lifting can present challenges in classical computation due to curse of dimensionality, it does not pose the same difficulties in quantum computing due to the use of qubits.

In recent work \cite{JLY2025fractionalPoisson}, we employ the Caffarelli-Silvestre extension technique \cite{Caffarelli2007Fractional} to numerically solve the fractional Poisson or Laplacian equation. This technique transforms the equation into a local elliptic PDE in one higher dimension, and we subsequently develop a new quantum algorithm to solve it.

In this paper, we focus on quantum computation for the time-fractional heat equation. Let $\Omega$ be an open, bounded and connected subset of $\mathbb{R}^d$ ($d\ge 1$) with Lipschitz boundary $\partial \Omega$. Given $\alpha\in (0,1)$, the fractional heat equation for $u$ is \cite{Khristenko23fractional}
\[\partial_t^{\alpha} u(t,x) - \Delta u(t,x) = 0,\]
where the fractional derivative is defined by
\[\partial_t^{\alpha} u = \frac{1}{\Gamma(1-\alpha)} \int_0^t (t-s)^{-\alpha} \partial_s u(s) \d s. \]

Due to the nonlocality and strong singularity in the problem,  classical numerical simulations of fractional PDEs encounter more  challenges compared to the local PDEs.
But on the other hand, as demonstrated in \cite{Eritz23timeFractional}, a time-fractional gradient flow can be  transformed into an integer-order one on an extended domain. This approach leads to a system of ODEs resulting from discretizations along both the extended and spatial domains.  Building on this extension technique, we develop a new quantum algorithm based on the Schr\"odingerization approach.

However, the original ODEs, which use the extended function at the collocation points as variables, have both positive and negative eigenvalues for the symmetric part of the coefficient matrix. This results in an unstable ODE system for the Schr\"odingerization approach, posing measurement difficulties. By applying a simple change of variables, we transform the original system into a stable one, which is well-suited for our quantum implementation and does not introduce additional complications when recovering the original variables for the fractional heat equation. A detailed complexity analysis shows that the quantum algorithm provides an {\it exponential} advantage over its classical counterparts, particularly in high-dimensional settings. Specifically, the quantum algorithm has a time complexity of $\widetilde{\mathcal{O}}(T^2d^4 h^{-8})$ for block-encoding input models, where $h$ is the spatial mesh size, $d$ is the number of dimensions, and $T$ is the evolution time.
In contrast, the classical method scales at least as $\mathcal{O}(N_t d h^{-(d+0.5)})$, where $N_t$ is the number of time steps, thereby demonstrating the exponential quantum advantage of our quantum algorithm with respect to the dimension $d$.

 The paper is structured as follows. In Section \ref{sec:fracGradFlows}, we provide a brief review of time-fractional gradient flows and introduce the extension technique discussed in \cite{Eritz23timeFractional}. Section \ref{sec:fracheat} applies this dimension extension to the fractional heat equation and derives the corresponding discrete ODE system. In Section \ref{sec:Schrfracheat}, we present the block encoding of the stiffness matrix on a hypercube domain, describe the process of extracting the solution vector, and analyze the query complexity of our Schr\"odingerization based approach. Numerical simulations and conclusions are presented in Section \ref{sec:numerical} and the final section, respectively.

\section{Fractional and extended gradient flows} \label{sec:fracGradFlows}

This section reviews a dimension-lifting technique introduced in \cite{Eritz23timeFractional} for the numerical solution of time-fractional gradient flows in Hilbert spaces. This method transforms a fractional flow into an equivalent integer-order gradient flow on an extended Hilbert space.

\subsection{Time-fractional gradient flows}

Let $T>0$ be a constant time and $X$ be a Banach space. We call a function $u: (0,T) \to X$ Bochner measurable if it can be approximated by a sequence of Banach-valued simple functions. The Bochner space $L^p(0, T; X)$ is defined as the equivalence class of Bochner measurable functions $u: (0,T) \to X$ such that $t \mapsto  \|u(t)\|_X^p$ is Lebesgue integrable. Similarly, we can define the Sobolev-Bochner space $W^{1,p}(0,T; X)$.

For $\alpha \in (0,1)$, we define the Riemann-Liouville integral operator of a function $u: (0,T) \to X$ as
\[\mathcal{I}_\alpha u(t) = \frac{1}{\Gamma(\alpha)} \int_0^t (t-s)^{-1+\alpha} u(s) \d s=: g_{\alpha} * u,\]
where $g_{\alpha} = t^{-1+\alpha}/\Gamma(\alpha)$ and $\Gamma(\alpha)$ is the Gamma function. The associated fractional derivative is defined by
\[\partial_t^{\alpha} u = \mathcal{I}_{1-\alpha} \partial_t u = g_{1-\alpha} * \partial_t u  =
\frac{1}{\Gamma(1-\alpha)} \int_0^t (t-s)^{-\alpha} \partial_s u(s) \d s. \]

Let $H$ be a Hilbert space. The time-fractional gradient flow in $H$ is defined as the following variational problem: Find $u \in L^p(0,T; X)$ such that
\begin{equation}\label{timefracflowVar}
(\partial_t^{\alpha} u, v)_H + \delta \mathcal{E}(u,v) = 0, \qquad v \in X
\end{equation}
for a given nonlinear energy functional $\mathcal{E}: X \to \mathbb{R}$, where $\delta \mathcal{E}: X \times X \to \mathbb{R}$ denotes its G\^{a}teaux derivative,
\[\delta \mathcal{E}(u,v) = \lim_{h \to 0} \frac{\mathcal{E}(u + h v) - \mathcal{E}(u)}{h}, \qquad u, v \in X.\]
The gradient of $\mathcal{E}$ is defined as $\nabla_H \mathcal{E}: X \to X'$ satisfying
\[\langle \nabla_H \mathcal{E}(u), v \rangle_{X'\times X} = \delta \mathcal{E}(u,v), \qquad v \in X.\]
With this notation, \eqref{timefracflowVar} can be written as
\[\partial_t^{\alpha} u = - \nabla_H \mathcal{E}(u) \quad \mbox{in}  \quad X', \qquad
\mbox{or} \qquad
\langle \partial_t^{\alpha} u + \nabla_H \mathcal{E}(u), v \rangle_{X'\times X} = 0, \qquad v\in  X.\]

 The variational problem is equipped with the initial data $ u_0 \in H $. As discussed in \cite{Eritz23timeFractional}, in general, one does not have $ u \in C([0,T]; H) $, and therefore cannot conclude that $ u(t) \to u_0 $ in $ H $ as $ t \to 0 $. However, the initial data are satisfied in the sense that $ g_{1-\alpha} * (u - u_0)(t) \to 0 $ in $ H $ as $ t \to 0 $, since $ g_{1-\alpha} * (u - u_0) \in C([0,T]; H) $. Furthermore, when $ \alpha > 1/p $, one has $ u - u_0 \in C([0,T]; X') $, so the initial condition is satisfied in the sense that $ u(t) \to u_0 $ in $ X' $. In fact,  the following well-posedness result of time-fractional gradient flows holds.

\begin{proposition}[Theorem 1 in \cite{Eritz23timeFractional}] \label{pro:fractionalgradientflow}
Let $X$ be a separable, reflexive Banach space that is compactly embedded in the separable Hilbert space $H$. Let $u_0\in H$, $\mathcal{E}\in C^1(X;\mathbb{R})$, and $\delta \mathcal{E}: X \times X \to \mathbb{R}$ be bounded and semicoercive in the sense
\begin{align*}
& \|\nabla_H\mathcal{E}(u)\|_{X'} \le C_0\|u\|_X^{p-1} + C_0, \\
& \delta\mathcal{E}(u, u) \ge C_1 \|u\|_X^p - C_2 \|u\|_H^2,
\end{align*}
for all $u \in X$, $p>1$, and $C_0, C_1, C_2$ are some positive constants. Moreover, assume that the realization of $\nabla_H\mathcal{E}$ as an operator from $L^p(0, T; X)$ to $L^{p'}(0,T;X')$ is weak-to-weak continuous. Then the time-fractional gradient flow $\partial_t^{\alpha} u = - \nabla_H \mathcal{E}(u)$ with $\alpha \in (0,1]$ admits a variational solution in the sense that
\[u \in L^p(0,T;X) \cap L^{(\infty)}(0,T;H) \cap W^{\alpha, p'}(0,T;X'), \qquad g_{1-\alpha} * (u-u_0) \in C([0,T]; H)\]
fulfills the variational form \eqref{timefracflowVar}.
\end{proposition}

We consider a simple case where the energy functional is taken as
\[\mathcal{E}(u) = \int_{\Omega} \frac12 |\nabla u(x)|^2 \d x,\]
where $X \subset H^1(\Omega)$ is regular enough. For the Hilbert space $H$, we list three choices considered in \cite{Eritz23timeFractional}:
\begin{itemize}
  \item $H = L^2(\Omega)$. The G\^{a}teaux derivative is
  \[\delta \mathcal{E}(u,v) = (-\Delta u, v), \qquad u, v \in X.\]

  \item $H = \dot{H}^1(\Omega)$, where
  \[\dot{H}^1(\Omega):= \Big\{ u \in H^1(\Omega): \int_{\Omega} u \d x = 0 \Big\},\]
  equipped with the scalar product $(\nabla \cdot, \nabla \cdot)_{L^2(\Omega)}$. The G\^{a}teaux derivative is
  \[\delta \mathcal{E}(u,v) = (u, v)_{\dot{H}^1(\Omega)} = (\nabla u, \nabla v)_{L^2(\Omega)}, \qquad u, v \in X.\]

  \item $H = \dot{H}^{-1}(\Omega)$, where $\dot{H}^{-1}(\Omega) = (\dot{H}^1(\Omega))'$ is the dual space of $\dot{H}^1(\Omega)$. The G\^{a}teaux derivative is
  \[\delta \mathcal{E}(u,v) = (\Delta^2 u,  v)_{\dot{H}^{-1}(\Omega)}, \qquad u, v \in X.\]
\end{itemize}
In particular, when $X\in C_c^\infty(\Omega)$, the respective gradient flows are
\begin{align}
& \partial_t^{\alpha} u = -u, \qquad H = H^1(\Omega), \nonumber \\
& \partial_t^{\alpha} u = \Delta u, \qquad H = L^2(\Omega),  \label{fracheat} \\
& \partial_t^{\alpha} u = -\Delta^2 u, \qquad H = \dot{H}^{-1}(\Omega).\nonumber
\end{align}

For the design of the quantum algorithm, we will focus on the second case, where the gradient of $ \mathcal{E} $ is given by
\begin{equation}\label{gradheat}
\nabla_H \mathcal{E}(u) = -\Delta u.
\end{equation}
which corresponds to the time-fractional heat equation.   As will be shown later, our discussion can be directly extended to other cases, as the discretization along the extended dimension remains the same. For other cases, we only need to adjust the spatial discretization to match the corresponding linear operator~-~such as the biharmonic operator in the third case.

\subsection{Augmented integer-order gradient flow}

 Let $ u: (0,T) \to H $ be a Bochner measurable function. Ref.~\cite{Eritz23timeFractional} introduced a dimension-lifting technique that transforms the time-fractional gradient flow into an integer-order gradient flow. To achieve this, the authors defined a function $ \tilde{u} $ on an extended domain $ \tilde{u}: (0,T) \times (0,1) \times \Omega \to \mathbb{R} $, interpreting it as a function mapping to a Hilbert space. Specifically,  define
\[
\tilde{u}: (0,T) \to L^2(0,1; H)
\]
such that $ \tilde{u}(t)(\theta, x) = \tilde{u}(t, \theta, x) $ for all $ t \in (0,T) $, $ \theta \in (0,1) $, and $ x \in \Omega $.

\begin{definition}
Let $\mathcal{E} :  X \to \mathbb{R}$ be a Fr\'echet differentiable energy functional that satisfies the assumptions in Proposition \ref{pro:fractionalgradientflow}. Consider the following two gradient flow problems:
\begin{itemize}
  \item the original time-fractional gradient flow in the Hilbert space $H$ with an initial $u_0 \in X$,
  \begin{equation}\label{originalflow}
  \partial_t^{\alpha} u = - \nabla_H \mathcal{E}(u), \qquad u(0) = u_0,
  \end{equation}
  \item and a higher-dimensional integer-order gradient flow in an augmented Hilbert space $\tilde{H}$ with zero initial,
  \begin{equation}\label{augmentedflow}
  \partial_t \tilde{u} = - \nabla_{\tilde{H}} \tilde{\mathcal{E}}(\tilde{u}), \qquad \tilde{u}(0) = 0.
  \end{equation}
\end{itemize}
Eq.  \eqref{augmentedflow} is called  the augmented gradient flow if one can find an appropriate form of the associated energy functional $\tilde{\mathcal{E}}$
and the Hilbert space $\tilde{H}$, depending on $\alpha$, which provide an equivalence of the variational solutions $u$ and $\tilde{u}$.
\end{definition}

To provide the form of the energy functional $\tilde{\mathcal{E}}$, introduce the following weighted Lebesgue spaces
\begin{align*}
& L_{\alpha}^2(0,1; H) = L^2(0,1; H; w_{\alpha}), \quad w_{\alpha} = g_{1-\alpha}(\theta) g_{\alpha}(1-\theta), \\
& L_{\alpha,0}^2(0,1; H) = L^2(0,1; H; w_{\alpha,0}), \quad w_{\alpha,0} = w_{\alpha}(\theta)(1-\theta), \\
& L_{\alpha,1}^2(0,1; H) = L^2(0,1; H; w_{\alpha,1}), \quad w_{\alpha,1} = w_{\alpha}(\theta) \theta
\end{align*}
for a given Hilbert space $H$. Here, the norm in $L^2(0,1; H; w)$ is defined by
\[(\tilde{u}, \tilde{v})_{L^2(0,1; H; w)} = \int_0^1(\tilde{u}(\theta), \tilde{v}(\theta) )_H  w(\theta) \d \theta. \]
Define
\[c_0(\theta) = \frac{1}{1-\theta},  \qquad c_1(\theta) = \frac{\theta}{1-\theta}.\]
Then it holds $w_{\alpha} = c_0w_{\alpha,0}$ and $w_{\alpha,1} = c_1 w_{\alpha,0}$. For any $\tilde{v} \in L_{\alpha}^2(0,1;H)$, define an integral operator $\mathcal{C}: L_{\alpha}^2(0,1;H) \to H$,
\[\mathcal{C} \tilde{v} = u_0 + (1,\tilde{v})_{L_{\alpha}^2(0,1)}
= u_0 + (c_0,\tilde{v})_{L_{\alpha,0}^2(0,1)} = u_0 + \int_0^1 w_{\alpha}(\theta) \tilde{v}(\theta) \d \theta,\]
and a functional $\mathcal{H}: L_{\alpha,1}^2(0,1;H) \to \mathbb{R}_{\ge 0}$,
\[\mathcal{H}(\tilde{v}) = \frac12 \|\tilde{v}\|_{L_{\alpha,1}^2(0,1;H)}^2
=\frac12 (c_1\tilde{v}, \tilde{v})_{L_{\alpha}^2(0,1; H)}
= \frac12\int_0^1 w_{\alpha,1}(\theta)\|\tilde{v}(\theta)\|_H^2 \d \theta.\]

The equivalence of the time-fractional gradient flow and its integer-order counterpart is described below.

\begin{proposition}[Theorem 2 in \cite{Eritz23timeFractional}] \label{pro:equivalence}
Let the assumptions from Proposition \ref{pro:fractionalgradientflow} hold.  Let $\tilde{X} = L_{\alpha}^2(0,1; X)$ and $\tilde{H}= L_{\alpha,0}^2(0,1; H)$ with $\alpha \in (0,1)$. The augmented energy functional is defined as
\[\tilde{\mathcal{E}} = \mathcal{E} \circ \mathcal{C} + \mathcal{H}.\]
Then, for any solution $u$ to the variational form of \eqref{originalflow} there exists a solution $\tilde{u}$ to the variational form of \eqref{augmentedflow}, and vice versa, such that the following equivalence of solutions holds:
\begin{equation}\label{uutilde}
u(t) = \mathcal{C} \tilde{u}(t) =  u_0 + \int_0^1 w_{\alpha}(\theta) \tilde{u}(t,\theta) \d \theta,
\end{equation}
\begin{equation}\label{utildeu}
\tilde{u}(t,\theta) = \mathcal{G}\{u\}(t,\theta):=
-c_0(\theta) \int_0^t \e^{-c_1(\theta)(t-s) }\nabla_H\mathcal{E}(u(s)) \d s.
\end{equation}
\end{proposition}

From \eqref{utildeu}, the augmented gradient flow \eqref{augmentedflow} can be expressed as
\begin{equation}\label{explicitaugment}
(1-\theta) \partial_t \tilde{u}(t,\theta) + \theta \tilde{u}(t,\theta) + \nabla_H \mathcal{E}(u(t)) = 0,
\end{equation}
where $ \nabla_H \mathcal{E}(u(t)) $ is given by \eqref{gradheat}, and $ u(t) $ can be written in terms of $ \tilde{u}(t) $ by substituting \eqref{uutilde} into the above equation.

\subsection{Rational approximation for the augmented gradient flow}

According to Proposition \ref{pro:equivalence}, one can derive an approximation of $u$ by solving \eqref{explicitaugment} at the quadrature points for the integral in \eqref{uutilde}. Let $\theta_k, \omega_k$ be the quadrature points and weights for \eqref{uutilde}. The solution to the fractional gradient flow is approximated by
\[u(t) \approx  u_0 + \sum_{k=1}^M \omega_k \tilde{u}^{(k)}(t),\]
where $\tilde{u}^{(k)}(t) = \tilde{u}(t,\theta_k)$ satisfies
\begin{equation}\label{naiveuk}
(1-\theta_k) \partial_t \tilde{u}^{(k)}(t) + \theta_k \tilde{u}^{(k)}(t) + \nabla_H \mathcal{E}(u(t)) = 0, \qquad k = 1,\cdots, M.
\end{equation}
However, in the weight function
\[w_{\alpha}(\theta) = \frac{1}{\Gamma(1-\alpha) \Gamma(\alpha)} \Big(\frac{\theta}{1-\theta}\Big)^{-\alpha} \frac{1}{1-\theta}, \]
the factor $\frac{\theta}{1-\theta}$ ranges from zero to infinity, which leads to a stiff problem for the numerical integration as well as the solution to \eqref{naiveuk}.

Using the change of variable $\lambda = \frac{\theta}{1-\theta}$, Eq.~\eqref{uutilde} can be written as
\[u(t) = u_0 + \frac{\sin (\pi \alpha) }{\pi} \int_0^{\infty} \frac{\lambda^{-\alpha}}{1+\lambda} \tilde{u}(t, \theta(\lambda)) \d \lambda.\]
According to the discussion in \cite{Eritz23timeFractional}, one can use the following quadrature rule for the above integral:
\begin{equation}\label{uexpand}
u(t) \approx  u_0 + \sum_{k=1}^M \frac{\omega_k}{1+\lambda_k} \tilde{u}^{(k)}(t) + \omega_{\infty} \tilde{u}^{(\infty)}(t),
\end{equation}
where $\theta_k = \frac{\lambda_k}{1+\lambda_k}$ and $\tilde{u}^{(k)}(t) = \tilde{u}(t, \theta(\lambda_k))$. The weights $\omega_k$ and the nodes $\lambda_k$ are given by the residues and the poles, respectively, in the rational approximation of the Laplace spectrum of the fractional kernel $g_{\alpha}$. In particular, one can use the adaptive Antoulas-Anderson (AAA) algorithm to determine the rational approximation of $\lambda^{-\alpha}$:
\begin{equation}\label{lambdafrac}
\lambda^{-\alpha} \approx \sum_{k=1}^M \frac{\omega_k}{\lambda + \lambda_k} + \omega_{\infty} ,  \qquad \lambda \in \Big[\frac{1}{T}, \frac{1}{\tau}\Big],
\end{equation}
where $0<\tau\ll T$ is a constant and $\omega_k, \lambda_k \ge 0$ (see \cite{Khristenko23fractional}).


The solution $u$ is approximated by a linear combination of $m+1$ modes $\tilde{u}^{(k)}(t) = \tilde{u}(t, \theta_k)$,
 and $\tilde{u}^{(k)}(t)$ satisfies the following system:
\begin{equation}\label{systemorig}
\begin{cases}
 \frac{1}{1+\lambda_k} \partial_t \tilde{u}^{(k)}(t) + \frac{\lambda_k}{1+\lambda_k}\tilde{u}^{(k)}(t) + \nabla_H \mathcal{E}(u(t)) = 0, \qquad k = 1,\cdots, M, \\
 \tilde{u}^{(\infty)}(t) + \nabla_H \mathcal{E}(u(t)) = 0.
\end{cases}
\end{equation}

\section{Dimension lifting for the time-fractional heat equation} \label{sec:fracheat}

This section applies the extension method from the previous section to solve the augmented system \eqref{systemorig}, with the gradient $\nabla_H \mathcal{E}(u(t))$ given by \eqref{gradheat}. This corresponds to the fractional heat equation
\begin{equation}\label{fracheat}
\partial_t^{\alpha} u(t,x) - \Delta u(t,x) = 0.
\end{equation}
We impose the Dirichlet boundary conditions with $x\in \Omega = (0,1)^d$, where $d$ is the spatial dimension.

\subsection{The vector form of the augmented system}

When $\nabla_H \mathcal{E}(u(t))$ is given by \eqref{gradheat}, we have
\begin{align}
& \frac{1}{1+\lambda_k} \partial_t \tilde{u}^{(k)}(t) + \frac{\lambda_k}{1+\lambda_k}\tilde{u}^{(k)}(t) -\Delta u(t) = 0, \qquad k = 1,\cdots, M,  \label{uk}\\
& \tilde{u}^{(\infty)}(t) -\Delta u(t) = 0, \label{uinf}
\end{align}
where $u(t)$ is approximated by \eqref{uexpand} and $\tilde{u}^{(k)}(0) = 0$ for $k=1,\cdots, m$.  For simplicity, we continue to use $ u(t) $ to denote the approximate solution and express \eqref{uexpand} as
\begin{equation}\label{uexpandequal}
u(t) =  u_0 + \sum_{i=1}^M \frac{\omega_i}{1+\lambda_i} \tilde{u}^{(i)}(t) + \omega_{\infty} \tilde{u}^{(\infty)}(t).
\end{equation}
The second equality \eqref{uinf} in the linear system along with \eqref{uexpandequal} gives
\[u(t) = ( I - \omega_{\infty} \Delta)^{-1} \Big( u_0 + \sum_{i=1}^M \frac{\omega_i}{1+\lambda_i} \tilde{u}^{(i)}(t) \Big) .\]
Therefore, we can derive an ODE system in terms of $\tilde{u}^{(k)}(t)$ ($k=1,\cdots,M)$:
\begin{equation}\label{systemu}
\begin{cases}
\partial_t \tilde{u}^{(k)}(t) + \lambda_k \tilde{u}^{(k)}(t)= (1+\lambda_k) \mathcal{L}_{\infty} \Big( u_0 + \sum_{i=1}^M \frac{\omega_i}{1+\lambda_i} \tilde{u}^{(i)}(t) \Big), \\
\tilde{u}^{(k)}(0) = 0,
\end{cases}
\end{equation}
where $\mathcal{L}_{\infty} = \Delta ( \mathcal{I} - \omega_{\infty} \Delta)^{-1}$ and $\mathcal{I}$ is an identity operator.

Introducing the vector notation
\[\tilde{\bb{u}}(t) = [\tilde{u}^{(1)}(t), \cdots, \tilde{u}^{(M)}(t)]^{\top},\]
where ``$\top$'' denotes the standard transpose, we can rewrite \eqref{systemu} in a compact form:
\begin{equation}\label{systemuvec}
\begin{cases}
\partial_t \tilde{\bb{u}}(t) = ( - \Lambda \mathcal{I} + C \mathcal{L}_{\infty})  \tilde{\bb{u}}(t) + \bb{f}, \\
\tilde{\bb{u}}(0) = \bb{0},
\end{cases}
\end{equation}
where $\Lambda = \text{diag}(\lambda_1,\cdots,\lambda_M)$, $C = (C_{ki})_{ M \times M}$ with $C_{ki} = (1+\lambda_k) \frac{\omega_i}{1+\lambda_i}$, and
\[\bb{f} = [ 1+\lambda_1, \cdots,  1+\lambda_n] ^{\top} \mathcal{L}_{\infty} u_0.\]

\subsection{The spatial discretization}

For simplicity, we only present the discretized version of \eqref{systemuvec} in one dimension. Here, we consider the use of the finite difference discretization with grid points given by
\[x_j = x_0 + j h, \qquad j = 0, 1,\cdots, N_x, \]
where $x_0 = 0 $ and $h = \frac{1}{N_x}$. At $x = x_j$ for $j = 1,\cdots, N_x-1$, the system \eqref{uk}-\eqref{uinf} can be discretized as
\begin{align*}
& \frac{1}{1+\lambda_k} \partial_t \tilde{u}_j^{(k)}(t) + \frac{\lambda_k}{1+\lambda_k}\tilde{u}_j^{(k)}(t) - \frac{u_{j-1}(t) - 2 u_j(t) + u_{j+1}(t)}{h^2} = 0, \qquad k = 1,\cdots, M,  \\
& \tilde{u}_j^{(\infty)}(t) -\frac{u_{j-1}(t) - 2 u_j(t) + u_{j+1}(t)}{h^2} = 0,
\end{align*}
where $\tilde{u}_j(t), u_j(t)$ are the approximate values of $\tilde{u}(t,x), u(t,x)$ at $x = x_j$.
Let $\tilde{\bb{u}}_h^{(k)}(t) = [\tilde{u}_1^{(k)}(t), \cdots, \tilde{u}_{N_x-1}^{(k)}(t)]^{\top}$ and $\bb{u}_h(t) = [u_1(t), \cdots, u_{N_x-1}(t)]^{\top}$. We have
\begin{align}
& \frac{1}{1+\lambda_k} \frac{\d }{\d t} \tilde{\bb{u}}_h^{(k)}(t) + \frac{\lambda_k}{1+\lambda_k}\tilde{\bb{u}}_h^{(k)}(t) = L_h \bb{u}_h(t) + \bb{b}_h(t), \qquad k = 1,\cdots, M,  \label{schemeFDM}\\
& \tilde{\bb{u}}_h^{(\infty)}(t) = L_h \bb{u}_h(t) + \bb{b}_h(t), \nonumber
\end{align}
where
\[L_h = \frac{1}{h^2}\begin{bmatrix}
-2       & 1      &       &    &    \\
 1       & -2     &  \ddots   &    & \\
         & \ddots & \ddots &  1  & \\
         &        & 1      & -2  &  1 \\
         &        &        &  1   & -2
\end{bmatrix}_{(N_x-1)\times (N_x-1)}, \qquad
\bb{b}_h(t) = \frac{1}{h^2}\begin{bmatrix}
u(t,0) \\
0 \\
\vdots \\
0 \\
u(t,1)
 \end{bmatrix}.\]
Accordingly, we can write the linear combination in \eqref{uexpand} as
\begin{equation}\label{uvecsum}
\bb{u}_h(t) =  \bb{u}_0 + \sum_{i=1}^M \frac{\omega_i}{1+\lambda_i} \tilde{\bb{u}}_h^{(i)}(t) + \omega_{\infty} \tilde{\bb{u}}_h^{(\infty)}(t).
\end{equation}
As in the operator form, by eliminating $\tilde{\bb{u}}_h^{(\infty)}(t)$, one gets
\begin{equation}\label{revoverut}
( I_h - \omega_{\infty} L_h)  \bb{u}_h(t) =  \bb{u}_0 + \sum_{i=1}^M \frac{\omega_i}{1+\lambda_i} \tilde{\bb{u}}_h^{(i)}(t) + \omega_{\infty} \bb{b}_h(t),
\end{equation}
where $I_h$ is the identity matrix of the same size as $L_h$.
Therefore, the spatial discretization of \eqref{systemuvec} is
\begin{align}
& \frac{\d }{\d t} \tilde{\bb{u}}_h^{(k)}(t) + \lambda_k \tilde{\bb{u}}_h^{(k)}(t) = (1+\lambda_k)   \sum_{i=1}^M \frac{\omega_i}{1+\lambda_i}  L_{\infty,h}  \tilde{\bb{u}}_h^{(i)}(t)   + (1+\lambda_k) \bb{f}_h(t),   \label{ukeq} \\
& \tilde{\bb{u}}_h^{(k)}(0) = \bb{0}, \qquad k = 1,\cdots, M, \label{ukeq0}
\end{align}
where
\[L_{\infty,h} = L_h ( I_h - \omega_{\infty} L_h)^{-1}, \qquad
\bb{f}_h(t) = \bb{b}_h(t) + L_{\infty,h} (\bb{u}_0+ \omega_{\infty} \bb{b}_h(t)). \]

Define
\[\tilde{\bb{U}}_h(t) = [\tilde{\bb{u}}_h^{(1)}(t); \cdots; \tilde{\bb{u}}_h^{(M)}(t)] ,\]
where ``;" indicates the straightening of $\{\tilde{\bb{u}}_h^{(i)}\}_{i\ge 1}$ into a column vector. We can express the ODE system as
\begin{equation}\label{systemuvecDiscretization}
\begin{cases}
\dfrac{\d }{\d t} \tilde{\bb{U}}_h(t) = ( - \Lambda \otimes I_h + C \otimes L_{\infty,h} ) \tilde{\bb{U}}_h(t) +  \tilde{\bb{F}}_h(t), \\
\tilde{\bb{U}}_h(0) = \bb{0},
\end{cases}
\end{equation}
where
\[\tilde{\bb{F}}_h = [ 1+\lambda_1, \cdots,  1+\lambda_M] ^{\top} \otimes \bb{f}_h(t).\]

In $d$ dimensions, the matrix $L_h$ should be replaced by
\begin{equation}\label{Lhd}
L_{h,d} = \underbrace{L_h\otimes I_h \otimes \cdots \otimes I_h}_{d~\text{matrices}} + I_h \otimes L_h\otimes \cdots \otimes I_h + \cdots + I_h \otimes I_h \otimes \cdots \otimes L_h.
\end{equation}
The boundary conditions can be deduced by following the same procedure in Section 4 of \cite{JLLY2024boundary}.

\begin{remark}
Compared to the standard heat equation, we need an additional $\log_2 m$ qubits for the fractional heat equation.
\end{remark}

\subsection{The modified ODE system}

Let $ \bb{a} = [1+\lambda_1, \cdots, 1+\lambda_M]^\top $ and $ \bb{b} = [\frac{\omega_1}{1+\lambda_1}, \cdots, \frac{\omega_M}{1+\lambda_M}]^\top $. Note the matrix $ C = \bb{a} \bb{b}^\top $ is a rank-one matrix. It is clear that the eigenvalues of $ C $ are $ c, 0, \dots, 0 $, where $ c = \bb{a}^\top \bb{b} = \omega_1 + \cdots + \omega_M > 0 $. However, this does not imply that the symmetric part
\[
D = \frac{C + C^\top}{2} = \frac{\bb{a} \bb{b}^\top + \bb{b} \bb{a}^\top}{2}
\]
is positive semi-definite. The matrix $ D $ is of rank 2, so it has two non-zero eigenvalues. Let $ \lambda $ be a non-zero eigenvalue of $ D $ with $ \bb{v} $ as the associated eigenvector. Then,
\[
D \bb{v} = \frac{1}{2} \left[\bb{a} (\bb{b}^\top \bb{v}) + \bb{b} (\bb{a}^\top \bb{v}) \right] = \lambda \bb{v}.
\]
This implies that $ \bb{v} $ is a linear combination of $ \bb{a} $ and $ \bb{b} $. Let $ \bb{v} = \alpha \bb{a} + \beta \bb{b} $. Substituting this into the above equation and performing simple manipulations, we obtain
\[
\begin{cases}
\frac{1}{2} \left[\alpha (\bb{b}^\top \bb{a}) + \beta (\bb{b}^\top \bb{b}) \right] = \lambda \alpha \\
\frac{1}{2} \left[\alpha (\bb{a}^\top \bb{a}) + \beta (\bb{a}^\top \bb{b}) \right] = \lambda \beta
\end{cases}
\quad \text{or} \quad
\begin{bmatrix}
\frac{\bb{b}^\top \bb{a}}{2} - \lambda   &  \frac{\bb{b}^\top \bb{b}}{2}  \\
 \frac{\bb{a}^\top \bb{a}}{2}    &  \frac{\bb{a}^\top \bb{b}}{2} - \lambda
\end{bmatrix}
\begin{bmatrix} \alpha \\ \beta \end{bmatrix} = \bb{0}.
\]
To get a non-zero solution, we require
\[
\text{det}\left(\begin{bmatrix}
\frac{\bb{b}^\top \bb{a}}{2} - \lambda   &  \frac{\bb{b}^\top \bb{b}}{2}  \\
 \frac{\bb{a}^\top \bb{a}}{2}    &  \frac{\bb{a}^\top \bb{b}}{2} - \lambda
\end{bmatrix}\right) = 0,
\]
which gives
\[
\lambda = \frac{\bb{a}^\top \bb{b} \pm \sqrt{(\bb{a}^\top \bb{a})(\bb{b}^\top \bb{b})}}{2}.
\]
Therefore, $ C $ has both a positive and a negative eigenvalue, which results in a unstable ODE system in \eqref{systemuvecDiscretization} for long-time evolution.

To overcome the issue, we denote $\hat{\bb{u}}_h^{(k)} = \frac{\omega_k}{1+\lambda_k}\tilde{\bb{u}}_h^{(k)}$ and rewrite \eqref{ukeq} as
\[\frac{\d }{\d t} \hat{\bb{u}}_h^{(k)}(t) + \lambda_k \hat{\bb{u}}_h^{(k)}(t) =   \omega_k \sum_{i=1}^M L_{\infty,h}  \hat{\bb{u}}_h^{(i)}(t) + \omega_k\bb{f}_h(t), \qquad k = 1,\cdots, M.\]
Define $\hat{\bb{U}}_h(t) = [\hat{\bb{u}}_h^{(1)}(t); \cdots; \hat{\bb{u}}_h^{(M)}(t)]$. We get
\begin{equation}\label{systemmodified}
\begin{cases}
\dfrac{\d }{\d t} \hat{\bb{U}}_h(t) = ( - \Lambda \otimes I_h + W \otimes L_{\infty,h} ) \hat{\bb{U}}_h(t) +  \hat{\bb{F}}_h(t), \\
\hat{\bb{U}}_h(0) = \bb{0},
\end{cases}
\end{equation}
where
\[W = \bb{\omega} \mathbbm{1} ^\top, \qquad \hat{\bb{F}}_h = \bb{\omega} \otimes \bb{f}_h(t), \qquad \bb{\omega} = [\omega_1,\cdots, \omega_M]^\top,\]
and $\mathbbm{1} \in \mathbb{R}^M$ is a column vector with all entries being 1. Let $D_{\omega} = \text{diag}(\bb{\omega})$ be a diagonal matrix. Noting that
\[D_{\omega}^{-1/2} W D_{\omega}^{1/2}  = D_{\omega}^{-1/2} \bb{\omega}\mathbbm{1}^\top D_{\omega}^{1/2}  = \bb{\omega}^{1/2} (\bb{\omega}^{1/2})^\top=:E_{\omega},\]
where $\bb{\omega}^{1/2} = [\omega_1^{1/2},\cdots, \omega_M^{1/2}]^\top$, we can rewrite \eqref{systemmodified} as
\begin{equation}\label{systemmodifiedfinal}
\begin{cases}
\dfrac{\d }{\d t} \bb{U}_h(t) = ( - \Lambda \otimes I_h + E_{\omega} \otimes L_{\infty,h} ) \bb{U}_h(t) +  \bb{F}_h(t), \\
\bb{U}_h(0) = \bb{0},
\end{cases}
\end{equation}
where
\[\bb{U}_h = (D_{\omega}^{-1/2} \otimes I_h) \hat{\bb{U}}_h, \qquad
\bb{F}_h = (D_{\omega}^{-1/2} \otimes I_h) \hat{\bb{F}}_h = \bb{\omega}^{1/2} \otimes \bb{f}_h(t).\]
This is a stable ODE system, as $ E_{\omega} $ is a symmetric rank-one matrix with
 \[\lambda = \sqrt{\bb{\omega}^\top \bb{\omega}} = \|\bb{\omega}\| >0 \]
 being the only non-zero eigenvalue and the symmetric part of $L_{\infty, h}$, i.e.,
 \[L_{\infty, h}^{\text{sym}} = \frac{1}{2}(L_{\infty, h} + L_{\infty, h}^\top)\]
 is negative definite.

For the stable system \eqref{systemmodifiedfinal}, if we let $\bb{U}_h = [\bb{u}_h^{(1)}; \cdots; \bb{u}_h^{(M)}]$, then the associated variables are
\[\bb{u}_h^{(k)} = \omega_k^{-1/2} \hat{\bb{u}}_h^{(k)} = \frac{\omega_k^{1/2}}{1+\lambda_k} \tilde{\bb{u}}_h^{(k)}, \qquad k = 1,\cdots, M.\]
From \eqref{uvecsum} we get
\begin{align}
( I_h - \omega_{\infty} L_h) \bb{u}_h(t)
& =  \bb{u}_0 + \sum_{i=1}^M \omega_i^{1/2} \bb{u}_h^{(i)}(t) + \omega_{\infty} \bb{b}_h(t) \nonumber \\
& = \bb{u}_0 + S \bb{U}_h + \omega_{\infty} \bb{b}_h(t), \label{utmodified}
\end{align}
where
\[S = [\omega_1^{1/2} I_h, \cdots, \omega_M^{1/2} I_h ]
 = (\bb{\omega}^{1/2})^\top \otimes I_h.\]

\subsection{Error estimate for the semi-discretization}

Let $\tilde{u}(t,x_j), u(t,x_j)$ be the exact solutions of $\tilde{u}(t,x), u(t,x)$ at $x = x_j$ for the augment system. We have
\begin{align*}
& \frac{1}{1+\lambda_k} \partial_t \tilde{u}^{(k)}(t,x_j) + \frac{\lambda_k}{1+\lambda_k}\tilde{u}^{(k)}(t,x_j)
 = \frac{u(t,x_{j-1}) - 2 u(t,x_j) + u(t,x_{j+1})}{h^2} + r(t,x_j),   \\
& \tilde{u}^{(\infty),x_j}(t) = \frac{u(t,x_{j-1}) - 2 u(t,x_j) + u(t,x_{j+1})}{h^2} + r(t,x_j),
\end{align*}
where
\[r(t,x_j) = - \frac{h^2}{12} u_{xxxx}(t,x_j) + \mathcal{O}(h^4).\]
Let $\tilde{\bb{u}}^{(k)}(t) = [\tilde{u}^{(k)}(t,x_1), \cdots, \tilde{u}^{(k)}(t,x_{N_x-1})]^{\top}$, $\bb{u}(t) = [u(t,x_1),
\cdots, u(t,x_{N_x-1})]^{\top}$ and
\[\bb{r}_h(t) = [r(t,x_1), \cdots, r(t,x_{N_x-1})]^{\top}.\]
Following the previous procedure,
one can derive the ODE system, as in \eqref{systemmodifiedfinal}, corresponding to the exact solutions of the augmented system:
\begin{equation}\label{systemmodifiedfinalexact}
\begin{cases}
\dfrac{\d }{\d t} \bb{U}(t) = ( - \Lambda \otimes I + E_{\omega} \otimes L_{\infty,h} ) \bb{U}(t) +  \bb{G}_h(t), \\
\bb{U}(0) = \bb{0},
\end{cases}
\end{equation}
where $\bb{G}_h(t)$ has the same form as $\bb{F}_h(t)$, but with $\bb{b}_h$ replaced by $\bb{b}_h + \bb{r}_h$.

 Throughout this article, we neglect the errors arising from the numerical integration in \eqref{uexpand}.

\begin{lemma}\label{lem:errODE}
Let $\Omega = (0,1)^d$. For every $t \in [0,T]$, assume that $u(t,\cdot) \in C^4(\overline{\Omega})$ and that the fourth-order derivative of $u(t,\cdot)$ is bounded in each direction. Then there holds
\[\frac{1}{\sqrt{N_d} } \|\bb{u}(t) - \bb{u}_h(t)\| \lesssim \frac{\|\bb{\omega}\|}{\|\bb{\lambda}\|_{\min}} d h^2,\]
where
$\|\bb{\lambda}\|_{\min} := \min \{\lambda_1,\cdots, \lambda_M\}$ and $N_d = (N_x-1)^d$ is the number of entries of $\bb{u}(t)$.
\end{lemma}
\begin{proof}
We only consider the one-dimension case.
Let $\bb{E}_h(t) = \bb{U}(t) - \bb{U}_h(t)$. It is simple to find that $\bb{E}_h(t)$ satisfies the following ODE system:
\[\begin{cases}
\dfrac{\d }{\d t} \bb{E}_h(t) = (- \Lambda \otimes I_h + E_{\omega} \otimes L_{\infty,h}) \bb{E}_h(t) +  \bb{g}_h(t), \\
\bb{E}_h(0) = \bb{0},
\end{cases}
\]
where $\bb{g}_h(t) = \bb{\omega}^{1/2} \otimes [(I + \omega_{\infty} L_{\infty,h}) \bb{r}_h(t)]$.
According to  Duhamel's principle,
\[\bb{E}_h(t) = \int_0^t \e^{A(t-s)} \bb{g}_h(s) \d s.\]

The matrix $L_h$ can be diagonalized as $L_h = PD P^\dag$, where $P$ is a unitary matrix and
\begin{equation}\label{DLh}
D= \text{diag}(d_1,\cdots,d_{N_x-1}), \qquad
d_j =  -\frac{4}{h^2} \sin^2 \frac{j}{2N_x} <0.
\end{equation}
Let $D_{\infty} = D(I-\omega_{\infty} D)^{-1}$. One has
\[(I\otimes P) \bb{E}_h(t) = \int_0^t  \e^{ (-\Lambda \otimes I_h + E_{\omega} \otimes D_{\infty}) (t-s)}  \Big(\bb{\omega}^{1/2} \otimes [(I_h + \omega_{\infty} D_{\infty}) P\bb{r}_h(t)] \Big) \d s.\]
This implies
\begin{align}
\|\bb{E}_h(t)\|
& \le \int_0^t  \e^{ - \|\bb{\lambda}\|_{\min}  (t-s)} \|\bb{\omega}\|^{1/2} \max_i \Big( 1 + \frac{\omega_{\infty} d_i}{1-\omega_{\infty} d_i}  \Big) \|\bb{r}_h(t)\|  \d s \nonumber \\
& \le \|\bb{\omega}\|^{1/2}  \int_0^t  \e^{ - \|\bb{\lambda}\|_{\min}  (t-s)}\|\bb{r}_h(t)\|  \d s
\lesssim \|\bb{\omega}\|^{1/2} \sqrt{N_x-1} h^2 \int_0^t  \e^{ - \|\bb{\lambda}\|_{\min}  (t-s)}   \d s \nonumber\\
& =  \|\bb{\omega}\|^{1/2} \sqrt{N_x-1} h^2  \frac{1- \e^{-\|\bb{\lambda}\|_{\min} t}}{\|\bb{\lambda}\|_{\min}}
\le \frac{\|\bb{\omega}\|^{1/2}}{\|\bb{\lambda}\|_{\min}} \sqrt{N_x-1} h^2, \label{ErrorE}
\end{align}
where we used that $u \in C^4([0,1])$ and $|u_{xxxx}|$ is bounded.

Define $\bb{e}_h = \bb{u} - \bb{u}_h$. From \eqref{utmodified} we obtain
\[
\| \bb{e}_h(t) \|
=   \| ( I - \omega_{\infty} L_h)^{-1} ((\bb{\omega}^{1/2})^\top \otimes I_h) \bb{E}_h  \| \le \|\bb{\omega}\|^{1/2} \|\bb{E}_h  \|,
\]
which along with \eqref{ErrorE} yields the desired estimate.
\end{proof}

\section{Schr\"odingerization for the time-fractional heat equation} \label{sec:Schrfracheat}

This section presents a Schr\"odingerization based quantum algorithm for solving the resulting ODE system \eqref{systemmodifiedfinal}. For brevity, we rewrite it using a more conventional notation:
\begin{equation}\label{ODElinear}
\begin{cases}
\dfrac{\d }{\d t} \bb{u}(t) = A \bb{u}(t) + \bb{b}(t), \\
\bb{u}(0) = \bb{u}_0,
\end{cases}
\end{equation}
where $\bb{u}(t) = \bb{U}(t)$, $\bb{u}_0 = \bb{0}$, and
\begin{equation}\label{coefficient}
A = - \Lambda \otimes I_h + E_{\omega} \otimes L_{\infty,h}, \qquad \bb{b}(t) = \bb{F}_h(t).
\end{equation}

\subsection{A review of the Schr\"odingerization method}

The detailed implementation of the Schr\"odingerization method is presented in \cite{JLMY2025SchOptimal}. Below, we provide a brief overview of how to transform the non-unitary ODE system \eqref{ODElinear} into a unitary dynamical system.

By introducing an auxiliary vector function $\bm{r}(t)$ that remains constant in time if $\bm{b}\neq 0$, the ODE system \eqref{ODElinear} can be rewritten as a homogeneous system
\begin{equation}\label{eq: homo Au}
	\frac{\d}{\d t} \bm{u}_f
	= A_f \bm{u}_f, \quad
	A_f  =  \begin{bmatrix}
		A &  B \\
		O &O
	\end{bmatrix},\quad
	\bm{u}_f(0) = \begin{bmatrix}
		\bm{u}_0 \\
		\bm{r}_0
	\end{bmatrix},
\end{equation}
where $B=\text{diag}\{b_0/\gamma_0, \cdots, b_{N-1}/\gamma_{N-1}\}$ and $\bm{r}_0 = \text{diag}\{\gamma_0, \cdots, \gamma_{N-1}\}$, with
\begin{equation}\label{eq:gamma}
	\gamma_i = T \sup_{t\in [0,T]} |b_i(t)|, \qquad i = 0,1,\cdots,N-1.
\end{equation}
Here, each $\sup_{t\in [0,T]} |b_i(t)|$ can be replaced by its upper bound and we set $b_i/\gamma_i = 0$ if $b_i(t) \equiv 0$.

For the Schr\"odingerization approach, we first decompose $A_f$ into a Hermitian term and an anti-Hermitian term:
\[A_f(t) = H_1(t) + \i H_2(t), \qquad \i = \sqrt{-1},\]
where
\[
H_1(t) = \frac{A_f(t)+A_f^{\dagger}(t)}{2} =  \begin{bmatrix}
	H_1^A & \frac{ B}{2}\\
	\frac{B^{\top}}{2} &O
\end{bmatrix}, \quad H_2(t) = \frac{A_f(t)-A_f^{\dagger}(t)}{2 \i} = \begin{bmatrix}
H_2^A & \frac{B}{2\i}\\
-\frac{ B^{\top}}{2\i} &O
\end{bmatrix},
\]
with
\[H_1^A = \frac{A+A^{\dagger}}{2} , \qquad H_2^A = \frac{A-A^{\dagger}}{2\i} .\]

Using the warped phase transformation $\bb{w}(t,p) = \e^{-p} \bb{u}_f(t)$ for $p\ge 0$ and symmetrically extending the initial data to $p<0$,  the ODE system \eqref{eq: homo Au} is  then transformed to a system of linear convection equations:
\begin{equation}\label{u2v}
\begin{cases}
 \dfrac{\d}{\d t} \bb{w}(t,p)  = - H_1(t) \partial_p \bb{w} + \i H_2(t) \bb{w}, \\
 \bb{w}(0,p) = \psi(p) \bb{u}_0,
 \end{cases}
\end{equation}
where $\psi(p)=\e^{-|p|}$. To get (near-) optimal cost, one can replaced it by a smooth extension (see \cite{JLMY2025SchOptimal}).
According to \cite[Theorem 3.1]{JLM24SchrInhom}, we can restore the solution $\bb{u}_f(t)$ by
\begin{equation}\label{eq:recovery}
\bm{u}_f = \e^p \bm{w}(t,p),\quad p\geq p^{\Diamond}   = \frac12.
\end{equation}

Suppose that the $p$-domain is $[-L,R]$, where $L,R>0$ are sufficiently large, partitioned by uniform grid points $-L = p_0<p_1<\cdots<p_{N_p} = R$.
Let the vector $\bb{w}_h$ be the collection of the function $\bb w$ at these grid points, defined more precisely as $
\bb{w}_h(t) = \sum_{k,i} \bb{w}_i(t,p_k) \ket{k,i}$, where $\bb{w}_i$ is the $i$-th entry of $\bb{w}$. Applying the discrete Fourier transform in the $p$ direction, one arrives at
\begin{equation}\label{heatww}
\frac{\d}{\d t} \bb{w}_h(t) = -\i (P_\mu \otimes  H_1 ) \bb{w}_h + \i (I\otimes H_2 ) \bb{w}_h , \quad
\bb{w}_h(0) = \bb{\psi} \otimes \bb{u}_0,
\end{equation}
where $\bb{\psi} = [\psi(p_0), \cdots, \psi(p_{N_p-1})]^{\top}$, $P_\mu$ is the matrix expression of the momentum operator $-\i\partial_p$, given by
\begin{equation}\label{Pmu}
P_\mu = \Phi D_\mu \Phi^{-1},  \qquad D_\mu = \text{diag}(\mu_0, \cdots, \mu_{N_p-1}).
\end{equation}
Here, $\mu_k = \frac{2\pi}{R+L} ( k - \frac{N_p}{2}) $ are the discrete Fourier modes and
\[\Phi = (\phi_{jl})_{N_p \times N_p} = (\phi_l(x_j))_{N_p\times N_p}, \qquad \phi_l(x) = \e ^{\i \mu_l (x+L)} \]
is the matrix representation of the discrete Fourier transform.
By a change of variables $\tilde{\bb{w}}_h = (\Phi^{-1} \otimes I)\bb{w}_h$, one derives the Hamiltonian system
\begin{equation}\label{generalSchr}
\frac{\d}{\d t} \tilde{\bb{w}}_h(t) = -\i H(t) \tilde{\bb{w}}_h(t) ,
\end{equation}
where $ H = D_\mu \otimes H_1 -  I \otimes H_2 $.

Both the near-optimal and optimal time complexity analysis of the Schr\"odingerization approach are presented in \cite{JLMY2025SchOptimal}, in which the following HAM-T oracle is assumed:
\begin{equation}\label{HAMT}
	(\bra{0^a} \otimes I) \text{HAM-T}_{H_\mu}(\ket{0^a} \otimes I)
	= \sum_{k=0}^{N_p-1} \ket{k}\bra{k} \otimes \frac{H_{\mu_k}}{\alpha_1 \mu_{\max} + \alpha_2},
\end{equation}
where $H_{\mu_k} = \mu_k H_1 - H_2$ and $\mu_{\max} = \max_k |\mu_k|$ represent the maximum absolute value among the discrete Fourier modes. This oracle only uses $\mathcal{O}(1)$ queries to block-encoding oracles for $H_1$ and $H_2$.

\begin{lemma} \label{lem:Schrcost}
Consider the linear ODE system \eqref{ODElinear}. Assume  $L$ and $R$ are large enough, $\triangle p = \mathcal{O}(\mu_{\max})$ is small enough and $\psi(p)$ in \eqref{u2v} is sufficiently smooth. Then there is a quantum algorithm that prepares an $\varepsilon$-approximation of the normalized solution $\ket{\bb{u}(T)}$ with $\Omega(1)$ probability and a flag indicating success, using
\[\widetilde{\mathcal{O}}\Big(\frac{\|\bm{u}(0)\|+T\|\bm{b}\|_{\text{smax}}}{\|\bb{u}(T)\|} \alpha_H  T \log^2 \frac{1}{\varepsilon}\Big)\]
queries to the HAM-T oracle in \eqref{HAMT} and
\[\mathcal{O}\Big(\frac{\|\bm{u}(0)\|+T\|\bm{b}\|_{\text{smax}}}{\|\bb{u}(T)\|}  \Big) \]
queries to the state preparation oracle for $\bb{\psi} \otimes \bb{u}_0$, where $\alpha_H \ge \max\{\alpha_1, \alpha_2\}$ and
\[\|\bm{b}\|_{\text{smax}} = \Big(\sum_{i=0}^{N-1} \Big(\sup_{t\in [0,T]} |b_i(t)| \Big)^2\Big)^{1/2}.\]
\end{lemma}

\subsection{Block-encoding of the coefficient matrix}

The matrices encountered are sparse, allowing us to construct all the involved matrices using sparse query models \cite{Gilyen2019QSVD, Chakraborty2019blockEncode, Lin2022Notes}. Given a matrix $A$, which is assumed to be sparse with at most $s_r$ nonzero entries in any row and at most $s_c$ nonzero entries in any column. In the following we set $s_r = s_c = s$. The sparse query model is described below.
\begin{definition}\label{def:sparsequery}
Let $A=(a_{ij})$ be an $n$-qubit matrix with at most $s$ non-zero elements in each row and column. Assume that $A$ can be accessed through the following oracles:
	\[
	O_r \ket{l} \ket{i} = \ket{r(i,l)} \ket{i}, \qquad O_c \ket{l} \ket{j} = \ket{c(j,l)} \ket{j},
	\]
	\[O_A \ket{0}\ket{i,j} = \Big( a_{ij} \ket{0} + \sqrt{1 - |a_{ij}|^2} \ket{1} \Big) \ket{i,j},\]
	where $r(i,l)$  and $c(j,l)$ give the $l$-th non-zero entry in the $i$-th row and $j$-th column.
\end{definition}
In the above definition, $O_A$ can be replaced by
\[\tilde{O}_A \ket{0^a} \ket{i,j} = \ket{\tilde{a}_{ij}} \ket{i,j},\]
where $\tilde{a}_{ij}$ is an $a$-bit binary representation of the $(i, j)$ matrix entry of $A$.

Compared to sparse encoding, block encoding is a more general input model for matrix operations on a quantum computer \cite{Gilyen2019QSVD,Chakraborty2019blockEncode,Lin2022Notes,ACL2023LCH2,JLMY2025SchOptimal}, which not only serves as an input model for quantum algorithms but also facilitates various matrix operations, enabling the block encoding of more complex matrices.

\begin{definition}\label{def:blockencoding}
	Suppose that $ A $ is an $ n $-qubit matrix and let $ \Pi = \bra{0^m} \otimes I_n$ with $ I_n$ being an $ n $-qubit identity matrix. If there exist positive numbers $ \alpha $ and $ \varepsilon $, as well as a unitary matrix $ U_A $ of $ (m+n) $-qubits, such that
	\[
	\| A - \alpha \Pi U_A \Pi^\dag \| = \| A - \alpha (\bra{0^m} \otimes I_n) U_A (\ket{0^m} \otimes I_n) \| \le \varepsilon,
	\]
	then $ U_A $ is called an $ (\alpha, m, \varepsilon) $-block-encoding of $ A $.
\end{definition}

\begin{figure}[!htb]
	\centering
	\includegraphics[scale=0.2]{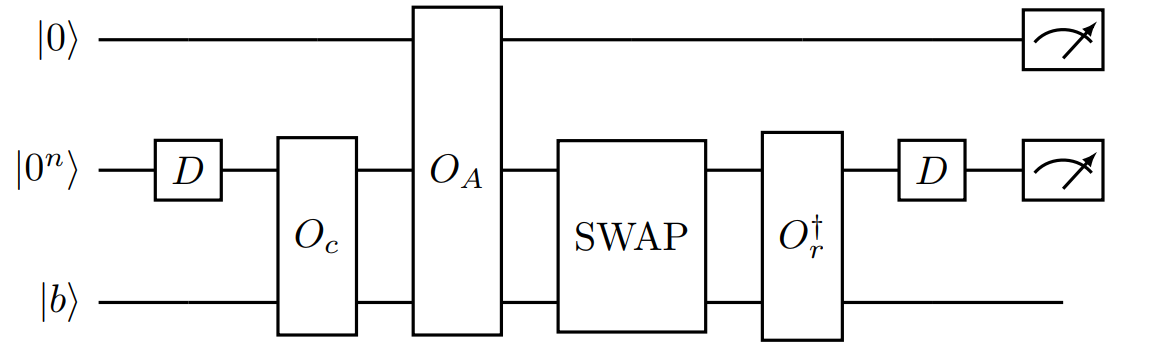}\\
	\caption{Quantum circuit for block encoding of sparse matrices.}\label{fig:BEASparseij}
\end{figure}

For a sparse matrix, we can construct its block encoding using the sparse query model \cite{Gilyen2019QSVD, Chakraborty2019blockEncode, Lin2022Notes, Lin2023timeMarching}.

\begin{lemma}[Block-encoding of sparse-access matrices] \label{lem:sparse2block}
	Let $A=(a_{ij})$ be an $n$-qubit matrix with at most $s$ non-zero elements in each row and column. Assume that  $\|A\|_{\max} = \max_{ij} |a_{ij}| \le 1$ and $A$ can be accessed through the sparse query model in Definition \ref{def:sparsequery}. Then we have an implementation of $(s, n+3, \epsilon)$-block-encoding of $A$ with a single use of $O_r$ and $O_c$, two uses of $\tilde{O}_A$ and additionally using $\mathcal{O}(n+\log^{2.5}(s/\epsilon))$ elementary gates and $\mathcal{O}(a+\log^{2.5}(s/\epsilon))$ ancilla qubits, where $a$ represents $a$-bit binary representations of the entries of $A$.
\end{lemma}

If $O_r$, $O_c$ and $O_A$ are exact, then we can implement an $(s, n+1, 0)$-block-encoding of $A$ with a single use of $O_r$, $O_c$ and $O_A$ and additionally using $\mathcal{O}(n)$ elementary gates and $\mathcal{O}(1)$ ancilla qubits. The circuit is shown in Fig.~\ref{fig:BEASparseij}, where the matrix $ D $ and $n$-qubit {\rm{SWAP}} operator are defined by
	\[
	D \ket{0^n} = \frac{1}{\sqrt{s}} \sum_{l \in [s]} \ket{l} = I_{n-\mathfrak{s}} \otimes \mathcal{H}^{\otimes \mathfrak{s}}, \quad s = 2^{\mathfrak{s}},\qquad
	 {\rm{SWAP}}\ket{i}\ket{j} = \ket{j}\ket{i}.
	\]
    Here $\mathcal{H}$ is a  Hadamard gate.

The following lemma describes how to perform matrix arithmetic on block-encoded matrices \cite{Gilyen2019QSVD,Chakraborty2019blockEncode}.

\begin{lemma}\label{lem:arithmeticBE}
	Let $A_i$ be $n$-qubit matrix for $i=1,2$. If $U_i$ is an $(\alpha_i, m_i, \varepsilon_i)$-block encoding of $A_i$ with gate complexity $T_i$, then
	\begin{enumerate}	
      \item $A_1 + A_2$ has an $(\alpha_1+\alpha_2, m_1+m_2, \alpha_1\varepsilon_2 + \alpha_2\varepsilon_1)$-block encoding that can be implemented with gate complexity $\mathcal{O}(T_1 + T_2)$.
	
      \item $A_1 A_2$ has an $(\alpha_1\alpha_2, m_1+m_2, \alpha_1\varepsilon_2 + \alpha_2\varepsilon_1)$-block encoding that can be implemented with gate complexity $\mathcal{O}(T_1 + T_2)$.

      \item $A_1^\dag$ has an $(\alpha_1, m_1, \varepsilon_1)$-block encoding that can be implemented with gate complexity $\mathcal{O}(T_1)$.

      \item $A_1\otimes A_2$ has an $(\alpha_1\alpha_2, m_1+m_2, \alpha_1\epsilon_2 + \alpha_2\epsilon_1)$-block encoding that can be implemented with gate complexity $\mathcal{O}(T_1+T_2)$.
	\end{enumerate}
\end{lemma}

We need the following result to invert a block encoded matrix.

\begin{lemma}\label{lem:inversecomplement}
Assume that $A$ can be accessed by its $(\alpha_A, n_A, \epsilon_A)$-block encoding $U_A$ with $\epsilon_A  \le \frac12\sigma_{\min}(A)$, where $\sigma_{\min}$ is the minimum singular value of $A$. Then one can construct a $(\frac{8}{3}\|A^{-1}\|, n_A+1, (1+2\|A^{-1}\|)\epsilon_A)$-block encoding of $A^{-1}$ using $\mathcal{O}(\alpha_A\|A^{-1}\|\log \frac{\|A^{-1}\|}{\epsilon_A})$ applications of (controlled-) $U_A$ and its inverse.
\end{lemma}
\begin{proof}
Ref.~\cite{Lin2023timeMarching} presented the result for $\epsilon_A = 0$ (see Section G.1 there). We extend the result to $\epsilon_A > 0$ using a similar argument.

For an odd function $f$ we define the singular value transformation $f^\Diamond$ as in \cite{Lin2023timeMarching}: if $A = W \Sigma V^\dag$, then $f^\Diamond(A): = W f(\Sigma) V^\dag$, where $W$ and $V$ are unitary matrices. We denote by $\tilde{A}$ the matrix that is exactly encoded by $U_A$, that is, $\alpha \Pi_1 U_A \Pi_1^\dag = \tilde{A}$, where $\Pi_1$ is the projection operator, $\|A - \tilde{A}\| \le \epsilon_A$ and $\|\tilde{A}\| \le \alpha_A$.

Let $f(x) = x^{-1}$. We have
\[(\tilde{A}/\alpha_A)^{-1} =  (f^\Diamond(\tilde{A}/\alpha_A) )^\dag.\]
Since $f(x)$ is not bounded by 1, we consider an odd polynomial $p(x)$ such that
\[\Big| p(x) - \frac{3\delta}{4x} \Big| \le \epsilon', \qquad x\in [-1, -\delta] \cup [\delta,1], \]
where $\delta \in (0,1)$ and $|p(x)| \le 1$ for all $x\in [-1,1]$. The existence of such an odd polynomial of degree $\mathcal{O}(\frac{1}{\delta} \log\frac{1}{\epsilon'})$ is guaranteed by using quantum singular value transformation (QSVT). Then we can implement $(p^\Diamond(\tilde{A}/\alpha_A) )^\dag$.

It is clear that
\begin{equation}\label{pdia}
\| (p^\Diamond(\tilde{A}/\alpha_A) )^\dag -\frac{3\delta}{4} (\tilde{A}/\alpha_A)^{-1} \| \le \epsilon'
\end{equation}
if all the singular values of $\tilde{A}/\alpha_A$ are in the interval $[\delta,1]$. Let $\tilde{A} = A+E$, where $\|E\| \le \epsilon_A$. Using the inequality
\[\|(I+B)^{-1}\| \le \frac{1}{1 - \|B\| }, \]
where $B$ satisfies $\|B\|<1$ and $I+B$ is invertible, we conclude that
\begin{align*}
\sigma_{\min}(\tilde{A})
& = \frac{1}{\|\tilde{A}^{-1}\|} = \frac{1}{\|(A+E)^{-1}\|}
 = \frac{1}{\|(I+A^{-1}E)^{-1}A^{-1}\|} \ge \frac{1 - \|A^{-1}\| \|E\| }{\|A^{-1}\|} \\
& \ge \frac{1 - \epsilon_A \|A^{-1}\| }{\|A^{-1}\|} = \frac12 \frac{1}{\|A^{-1}\|} = \frac12 \sigma_{\min}(A).
\end{align*}
This indicates that the singular values of $\tilde{A}/\alpha_A$ are contained in $[\frac{1}{2\|A^{-1}\|}, \|\tilde{A}\| ]$. Therefore, we can choose $\delta =\frac{1}{2\alpha_A \|A^{-1}\|}$.

The matrix $\tilde{A}$ has a $(1, n_A, 0)$-block encoding. Using QSVT, a $(1, n_A+1, 0)$-block encoding of $p^\Diamond(\tilde{A}/\alpha_A)$ can be constructed. We denote this block encoding by $\mathcal{U}$. Then by Eq.~\eqref{pdia},
\[\| \frac{4}{3\delta \alpha_A} (p^\Diamond(\tilde{A}/\alpha_A) )^\dag - \tilde{A}^{-1} \| \le \frac{4}{3\delta \alpha_A} \epsilon'.\]
According to the discussion in \cite{Lin2023timeMarching}, $\mathcal{U}^\dag$ is a $(4\|\tilde{A}^{-1}\|/3, n_A+1, \tilde{\epsilon})$-block encoding of $\tilde{A}^{-1}$, where
\[\tilde{\epsilon} = \frac{4}{3\delta \alpha_A} \epsilon' = \frac{8}{3 }\|A^{-1}\| \epsilon' .\]
The total number of queries to $U_A$ and its inverse is
\[\mathcal{O}\Big(\frac{1}{\delta} \log \frac{1}{\epsilon'} \Big)
= \mathcal{O}\Big(\alpha_A \|A^{-1}\| \log \frac{1}{\epsilon'}\Big).\]

It is simple to find that
\[\tilde{A}^{-1} - A^{-1} = (A+E)^{-1} - A^{-1} = (A+E)^{-1}( A - (A+E) ) A^{-1} = -(A+E)^{-1} E  A^{-1}, \]
yielding
\begin{align*}
\| \tilde{A}^{-1} - A^{-1}  \|
& \le \|(A+E)^{-1}\| \| E\| \|A^{-1}\| \le \epsilon_A \|(A+E)^{-1}\| \|A^{-1}\| \\
& \le \epsilon_A \frac{\|A^{-1}\|}{1-\|EA^{-1}\|} \|A^{-1}\| \le 2 \epsilon_A \|A^{-1}\|^2.
\end{align*}
The proof is finished by taking $\epsilon' = 3\epsilon_A/(8\|A^{-1}\|)$.
\end{proof}

Note that if $\|A^{-1}\|$ is unknown, it can be replaced by an upper bound of $\|A^{-1}\|$. Additionally, it is reasonable to assume that $\alpha_A = \mathcal{O}(\|A\|)$, implying that $\alpha_A \|A^{-1}\| = \mathcal{O}(\kappa_A)$, where $\kappa_A$ denotes the condition number of $A$.

With the help of the above preparations, we are in a position to block encode the coefficient matrix $ A $ in \eqref{coefficient}.
For simplicity, we assume $ M = 2^m $ and $ N_x - 1 = 2^{n_x} $, where $ M $ represents the number of equations in \eqref{naiveuk}.

According to Lemmas \ref{lem:sparse2block} and \ref{lem:arithmeticBE}, we can construct
\begin{itemize}
  \item an $(\alpha_{\Lambda}, m+3, \epsilon_1)$-block-encoding of $\Lambda$, where $\alpha_{\Lambda} = \|\bb{\lambda}\|_{\max}$, with gate complexity
      \[T_1 = \mathcal{O}(m+\log^{2.5}(1/\epsilon_1)).\]
  \item an $(\alpha_{L_h}, n_x+3, \epsilon_2)$-block-encoding of $L_h$, where $\alpha_{L_h} \simeq 3 h^{-2}$, with gate complexity
      \[T_2 = \mathcal{O}(n_x+\log^{2.5}(3/\epsilon_2)).\]
\end{itemize}
For the matrix $E_{\omega} = \bb{\omega}^{1/2} (\bb{\omega}^{1/2})^\top$, we can first construct an $(\alpha_{1/2}, m+3, \epsilon_3')$-block-encoding of $\bb{\omega}^{1/2}$, which leads to an $(\alpha_{E_{\omega}}, 2(m+3), \epsilon_3)$-block-encoding of $E_{\omega}$, where $\alpha_{E_{\omega}} = \alpha_{1/2}^2 \simeq \|\bb{\omega}\|$ and $\epsilon_3 = 2\alpha_{1/2}\epsilon_3'$, with gate complexity
 \[T_3 = \mathcal{O}(m+\log^{2.5}(1/\epsilon_3)).\]

\begin{theorem} \label{thm:inputmodel}
Given an $(\alpha_{\Lambda}, n_1, \epsilon_1)$-block encoding of $\Lambda$, an $(\alpha_{L_h}, n_2, \epsilon_2)$-block encoding of $L_h$ and an $(\alpha_{E_{\omega}}, n_3, \epsilon_3)$-block encoding of $E_{\omega}$, where
\begin{align*}
& n_1 = m+3, \qquad n_2 = n_x+3, \qquad n_3 = 2(m+3), \\
& \alpha_{\Lambda} \simeq \|\bb{\lambda}\|_{\max}, \qquad \alpha_{L_h} \simeq h^{-2}, \qquad \alpha_{E_{\omega}} \simeq \|\bb{\omega}\|,
\end{align*}
we can construct
\begin{itemize}
  \item an $(\alpha_{\text{inv}}, n_{\text{inv}}, \epsilon_{\text{inv}})$-block encoding, denoted by $\mathcal{U}_{\text{inv}}$, of $( I_h - \omega_{\infty} L_{h,d})^{-1}$, where
  \begin{align}
  & \alpha_{\text{inv}} = \mathcal{O} (1 + \omega_{\infty} d h^{-2} ),  \label{alphainverse}\\
  & n_{\text{inv}} = (d+d^2) n_2+1,  \nonumber \\
  & \epsilon_{\text{inv}} = \mathcal{O} (d h^{-2} (1 + \omega_{\infty} d h^{-2}) \epsilon_2), \nonumber
  \end{align}
  with gate complexity
  \[T_{\text{inv}} = \mathcal{O}\Big( dn_x (1 + \omega_{\infty} d h^{-2})^2\log^{3.5}\frac{1}{\epsilon_2} \Big);\]

  \item an $(\alpha_A, n_A, \epsilon_A)$-block encoding, denoted by $\mathcal{U}_A$, of the coefficient matrix $A = - \Lambda \otimes I_h + E_{\omega} \otimes L_{\infty,h}$
of \eqref{ODElinear} in $d$ dimensions, where
  \begin{align}
  & \alpha_A = \mathcal{O}(\|\bb{\lambda}\|_{\max} + \|\bb{\omega}\| dh^{-2}(1 + \omega_{\infty} d h^{-2})), \label{alphaA}\\
  & n_A = n_1 + (d+3d^2) n_2 + n_3 + 1,  \nonumber\\
  & \epsilon_A = \mathcal{O}\Big( dh^{-2}(1 + \omega_{\infty} d h^{-2})
  (\|\bb{\omega}\|\epsilon_1 + \|\bb{\lambda}\|_{\max}\|\bb{\omega}\|d h^{-2} \epsilon_2 +  \|\bb{\lambda}\|_{\max}\epsilon_3  )  \Big), \nonumber
  \end{align}
  with gate complexity
  \[T_A = \mathcal{O}(m + \log^{2.5}\frac{1}{\epsilon_1} +\log^{2.5}\frac{1}{\epsilon_3} + dn_x (1 + \omega_{\infty} d h^{-2})^2\log^{3.5}\frac{1}{\epsilon_2}).\]
\end{itemize}
\end{theorem}
\begin{proof}
In $d$ dimensions, $L_{\infty,h} = L_{h,d} ( I_h - \omega_{\infty} L_{h,d})^{-1}$, where $L_{h,d}$ is defined in \eqref{Lhd}.

We bound the complexity for $( I_h - \omega_{\infty} L_{h,d})^{-1}$ in the following steps:
\begin{enumerate}
  \item It is obvious that each term in the summation of \eqref{Lhd} has an $(\alpha_{L_h}, d n_2, \epsilon_2)$-block encoding with gate complexity $T_2$. By Lemma \ref{lem:arithmeticBE}, $L_{h,d}$ has a $(d \alpha_{L_h}, d^2 n_2, d \alpha_{L_h} \epsilon_2)$-block encoding with gate complexity $dT_2$.

  \item Using Lemmas \ref{lem:sparse2block} and \ref{lem:arithmeticBE}, we can construct a $(1 + \omega_{\infty} d \alpha_{L_h}, dn_2+d^2 n_2, d \alpha_{L_h} \epsilon_2)$-block encoding $U_{h,d}$ of $I_h - \omega_{\infty} L_{h,d}$ since $I_h$ has a $(1, dn_2, 0)$-block encoding, with gate complexity $\mathcal{O}(dT_2)$.

  \item Using Lemma \ref{lem:inversecomplement}, we have an $(\alpha_{\text{inv}}, n_{\text{inv}}, \epsilon_{\text{inv}})$-block encoding of $( I_h - \omega_{\infty} L_{h,d})^{-1}$, where
  \begin{align*}
  & \alpha_{\text{inv}} = \|(I_h - \omega_{\infty} L_{h,d})^{-1}\|  = \mathcal{O} (1 + \omega_{\infty} d h^{-2} ), \\
  & n_{\text{inv}} = (d+d^2)n_2+1, \\
  & \epsilon_{\text{inv}} = (3 + 2\omega_{\infty} d h^{-2}) d h^{-2} \epsilon_2,
  \end{align*}
  with
  \begin{align*}
  & \mathcal{O}\Big( (1 + \omega_{\infty} d \alpha_{L_h})\|(I_h - \omega_{\infty} L_{h,d})^{-1}\|\log
  \frac{\|I_h - \omega_{\infty} L_{h,d}\|^{-1}}{d \alpha_{L_h} \epsilon_2}  \Big) \\
  & = \mathcal{O}\Big( (1 + \omega_{\infty} d h^{-2})^2\log
  \frac{1 + \omega_{\infty} d h^{-2}}{d h^{-2} \epsilon_2}  \Big)
  = \mathcal{O}\Big( (1 + \omega_{\infty} d h^{-2})^2\log
  \frac{1}{\epsilon_2}  \Big)
  \end{align*}
  applications of $U_{h,d}$ and its inverse. This implies that the gate complexity is
  \[\mathcal{O}\Big( (1 + \omega_{\infty} d h^{-2})^2dT_2  \log
  \frac{1}{\epsilon_2} \Big).\]
\end{enumerate}

For the complexity for $A = - \Lambda \otimes I_h + E_{\omega} \otimes L_{\infty,h}$, we have the following steps:
\begin{enumerate}[(a)]
  \item According to the previous discussion, one can construct an $(\alpha_{\infty,h}, n_{\infty,h}, \epsilon_{\infty,h})$-block encoding of $L_{\infty,h} = L_{h,d} ( I_h - \omega_{\infty} L_{h,d})^{-1}$, where
    \begin{align*}
  & \alpha_{\infty,h} = d \alpha_{L_h} \alpha_{\text{inv}} = \mathcal{O} (dh^{-2}(1 + \omega_{\infty} d h^{-2})  ), \\
  & n_{\infty,h} = d^2 n_2 + n_{\text{inv}} =  (d+2d^2) n_2+1, \\
  & \epsilon_{\infty,h} = \alpha_{\text{inv}} d \alpha_{L_h} \epsilon_2 + d \alpha_{L_h} \epsilon_{\text{inv}} =
  \mathcal{O} ( d^2 h^{-4} (1 + \omega_{\infty} d h^{-2})  \epsilon_2),
  \end{align*}
  with gate complexity
  \[T_{\infty,h} = T_{\text{inv}} + \mathcal{O}(dT_2) = \mathcal{O}( T_{\text{inv}} ).\]

  \item The matrix $E_{\omega} \otimes L_{\infty,h}$ has an
  \[(\alpha_{E_{\omega}}\alpha_{\infty,h}, n_3+n_{\infty,h}, \alpha_{E_{\omega}}\epsilon_{\infty,h} + \alpha_{\infty,h}\epsilon_3)\]
  block encoding with gate complexity
  \[T_3 + T_{\infty,h} = \mathcal{O}(m+\log^{2.5}\frac{1}{\epsilon_3} + dn_x (1 + \omega_{\infty} d h^{-2})^2\log^{3.5}\frac{1}{\epsilon_2}).\]

  \item The matrix $- \Lambda \otimes I_h$ has an $(\alpha_{\Lambda}, n_1 + dn_2, \epsilon_1)$-block encoding, implying that $A = - \Lambda \otimes I_h + E_{\omega} \otimes L_{\infty,h}$ has an $(\alpha_A, n_A, \epsilon_A)$-block encoding, where
      \begin{align*}
  & \alpha_A = \alpha_{\Lambda} + \alpha_{E_{\omega}}\alpha_{\infty,h}
  = \mathcal{O}(\|\bb{\lambda}\|_{\max} + \|\bb{\omega}\| dh^{-2}(1 + \omega_{\infty} d h^{-2})), \\
  & n_A = n_1 + dn_2 + n_3+ n_{\infty,h} =  n_1 + (d+3d^2) n_2 + n_3 + 1, \\
  & \epsilon_A = \alpha_{\Lambda} (\alpha_{E_{\omega}}\epsilon_{\infty,h} + \alpha_{\infty,h}\epsilon_3) + \alpha_{E_{\omega}}\alpha_{\infty,h}\epsilon_1 \\
  & \quad = \mathcal{O}\Big( dh^{-2}(1 + \omega_{\infty} d h^{-2})
  (\|\bb{\omega}\|\epsilon_1 + \|\bb{\lambda}\|_{\max}\|\bb{\omega}\|d h^{-2} \epsilon_2 +  \|\bb{\lambda}\|_{\max}\epsilon_3  )  \Big).
  \end{align*}
\end{enumerate}
The gate complexity is
\[T_A = T_1 + T_3 + T_{\infty,h} = \mathcal{O}(m + \log^{2.5}\frac{1}{\epsilon_1} +\log^{2.5}\frac{1}{\epsilon_3} + dn_x (1 + \omega_{\infty} d h^{-2})^2\log^{3.5}\frac{1}{\epsilon_2}). \]
This completes the proof.
\end{proof}

\begin{remark}
When $\omega_\infty \neq 0$, the gate complexity for the input model is $\widetilde{\mathcal{O}}(d^3 h^{-4})$ in terms of the dimension $d$ and the inverse of $h$, which is comparable to that of the classical input model, where $\widetilde{\mathcal{O}}$ represents an upper bound that omits logarithmic factors.
\end{remark}

It is evident that the HAM-T oracle in \eqref{HAMT} can be constructed with  $\mathcal{O}(1)$ queries to the oracle $U_A$, with
\[\alpha_1, \alpha_2 \simeq \alpha_A + \frac{T}{2} =: \alpha_H.\]
For simplicity, we assess the complexity based on the queries to $\mathcal{U}_A$ in this article.

\subsection{Recover the solution from the augmented system} \label{subsec:recover}

For simplicity, we assume homogeneous boundary conditions in the following discussion.
The approximate solution $\bb{u}_h(T)$ is given by
\begin{equation}\label{uh}
(I_h -\omega_{\infty}L_{\infty,h}) \bb{u}_h(T) = \bb{u}_0 + \sum_{i=1}^M \omega_i^{1/2} \bb{u}^{(i)}_h(T) = \bb{u}_0 + S\bb{U}_h(T)=:\bb{u}_0 + \bb{u}_T,
\end{equation}
where $S = (\bb{\omega}^{1/2})^{\top} \otimes I_h$.

One can apply the LCU procedure to combine the two terms on the right-hand side of \eqref{uh}. To this end, we introduce the following unitaries:
\begin{itemize}
  \item For $\bb{u}_0$, we assume the state preparation oracle:
  \[\ket{0^{n_a}} \ket{0^{ m}} \ket{0^{d n_x}} \quad \xrightarrow{ I^{\otimes (n_a+m) } \otimes O_{\text{prep}} } \quad  \frac{1}{\eta_0} \ket{0^{n_a+m}} \otimes \bb{u}_0,\]
  where $\eta_0 = \|\bb{u}_0\|$ and $n_a$ is the number of ancilla qubits. For brevity, we denote $V_0$ to be the unitary.

  \item For $\bb{U}_h(T)$, applying the Schr\"odingerization approach, we can get
  \begin{equation}\label{Vunitary}
   \ket{0^{dn_x+m}} \quad \xrightarrow{ V } \quad  \ket{\bb{U}_h(T)},
  \end{equation}
  where the potential involvement of ancilla qubits has been neglected.

  \item For $S\bb{U}_h(T)$, assume the block encoding $U_S$ of $S$ with the encoding constant $ \alpha_S \simeq \|\bb{\omega}\|^{1/2} $. Then,
  \begin{align*}
  \ket{0^{n_a}} \ket{0^{d n_x+ m}}
  &  \quad \xrightarrow{ I^{\otimes n_a} \otimes V } \quad  \ket{0^{n_a}}\ket{\bb{U}_h(T)} \\
  &  \quad \xrightarrow{ ~~~U_S~~~ } \quad \frac{1}{\alpha_S} \ket{0^{n_a}}S\ket{\bb{U}_h(T)} + \ket{\bot} \\
  & \quad = \frac{1}{\eta_1 } \ket{0^{n_a+m}}S\bb{U}_h(T) + \ket{\bot},
   \end{align*}
  where $\eta_1 = \alpha_S \|\bb{U}_h(T)\|$.  We denote $V_1$ to be the unitary and let $\bb{u}_T = S\bb{U}_h(T)$.
\end{itemize}

Let $R_t$ be a single-qubit rotation such that
$R_t\ket{0} = \frac{1}{\sqrt{\eta_0 + \eta_1}}( \sqrt{\eta_0} \ket{0} + \sqrt{\eta_1} \ket{1})$.
Let $n_s = dn_x + m$. The desired LCU procedure is given as follows:
    \begin{align}
\ket{0} \otimes \ket{0^{n_a}} \otimes \ket{0^{n_s}}
&  \quad \xrightarrow{ R_t \otimes I^{\otimes n_a} \otimes I^{\otimes n_s} }
\quad
    \frac{1}{\sqrt{\eta_0 + \eta_1}}( \sqrt{\eta_0} \ket{0} + \sqrt{\eta_1} \ket{1}) \otimes \ket{0^{n_a}} \otimes \ket{0^{n_s}} \nonumber\\
&  \quad \xrightarrow{\ketbra{0}{0}\otimes V_0 + \ketbra{1}{1}\otimes V_1}
\quad
   \frac{1}{\sqrt{\eta_0} \sqrt{\eta_0 + \eta_1}}\ket{0} \otimes \ket{0^{n_a+m}} \otimes \bb{u}_0 \nonumber\\
&  \hspace{4cm} +  \frac{1}{\sqrt{\eta_1} \sqrt{\eta_0 + \eta_1}}\ket{1} \otimes \ket{0^{n_a}} \otimes \bb{u}_T + \ket{\bot}\nonumber\\
& \quad \xrightarrow{~R_t^\dag \otimes I^{\otimes n_a} \otimes I^{\otimes n_s}~}
\quad
\frac{1}{\eta_0 + \eta_1}\ket{0}\otimes \ket{0^{n_a}}
\otimes ( \bb{u}_0 + \bb{u}_T) + \ket{\bot}  , \label{u0plusuT}
\end{align}
where
\[\eta_0 + \eta_1 = \|\bb{u}_0\| + \alpha_S \|\bb{U}_h(T)\| \simeq \|\bb{u}_0\| + \|\bb{\omega}\|^{1/2} \|\bb{U}_h(T)\|.\]

The explicit time complexity necessitates an estimation of $\|\bb{U}_h(T)\|$.
\begin{lemma}\label{lem:boundUh}
Suppose that the time-fractional heat equation \eqref{fracheat} is imposed with homogeneous Dirichlet boundary conditions. Then the solution to \eqref{systemmodifiedfinal} satisfies
\[\frac{\|\bb{\omega}\|^{1/2}\|\bb{u}_0\|}{\|\bb{\lambda}\|_{\max} + \|\bb{\omega}\|/\omega_{\infty}} \lesssim \|\bb{U}_h(T)\| \lesssim \frac{\|\bb{\omega}\|^{1/2}\|\bb{u}_0\|}{\|\bb{\lambda}\|_{\min}}.\]
\end{lemma}
\begin{proof}
By Duhamel's principle, the solution to \eqref{systemmodifiedfinal} can be expressed as
\[\bb{U}_h(t) = \int_0^t \e^{A(t-s)} \d s\bb{F}_h ,\]
where
\[A = - \Lambda \otimes I  + E_{\omega} \otimes L_{\infty,h}, \qquad \bb{F}_h = \bb{\omega}^{1/2} \otimes ( L_{\infty,h} \bb{u}_0). \]
As in the proof of Lemma \ref{lem:errODE},
\[(I \otimes P) \bb{U}_h(t) = \int_0^t  \e^{ (-\Lambda \otimes I  + E_{\omega} \otimes D_{\infty}) (t-s)}  \d s \Big(\bb{\omega}^{1/2} \otimes (D_{\infty} P\bb{u}_0) \Big),\]
where $D_{\infty} = D(I -\omega_{\infty} D)^{-1}$ with $D$ defined in \eqref{DLh}. Noting that $d_i<0$, one obtains
\begin{align*}
\|\bb{U}_h(t)\|
& \le \int_0^t  \e^{ - \|\bb{\lambda}\|_{\min}  (t-s)} \|\bb{\omega}\|^{1/2} \max_i \Big|\frac{d_i}{1-\omega_{\infty} d_i}\Big| \|\bb{u}_0\|  \d s \\
& \lesssim \|\bb{\omega}\|^{1/2}\|\bb{u}_0\| h^{-2} \int_0^t  \e^{ - \|\bb{\lambda}\|_{\min}  (t-s)}   \d s \\
& = \|\bb{\omega}\|^{1/2}\|\bb{u}_0\|  \frac{1- \e^{-\|\bb{\lambda}\|_{\min} t}}{\|\bb{\lambda}\|_{\min}}
\le \frac{\|\bb{\omega}\|^{1/2}}{\|\bb{\lambda}\|_{\min}} \|\bb{u}_0\|.
\end{align*}

On the other hand,
\[\Big(\int_0^t  \e^{ A (t-s)}  \d s\Big)^{-1}  \bb{U}_h(t) =  \bb{\omega}^{1/2} \otimes (L_{\infty,h} \bb{u}_0),\]
which gives
\begin{align*}
\|\bb{\omega}\|^{1/2} \| L_{\infty,h} \bb{u}_0\|
 = \|((\e^{At}-I)A^{-1})^{-1}\bb{U}_h(t)\|  \le \|((\e^{At}-I)A^{-1})^{-1}\| \|\bb{U}_h(t)\|.
\end{align*}
This yields
\begin{align*}
\|\bb{U}_h(t)\|
& \ge \|\bb{\omega}\|^{1/2} \| L_{\infty,h} \bb{u}_0\|\|((\e^{At}-I)A^{-1})^{-1}\|^{-1} \\
& = \|\bb{\omega}\|^{1/2} \| L_{\infty,h} \bb{u}_0\| \sigma_{\min}((\e^{At}-I)A^{-1}) \\
& \ge \|\bb{\omega}\|^{1/2} \sigma_{\min}( L_{\infty,h} ) \| \bb{u}_0\| \min_i \Big|\frac{1-\e^{\mu_i t}}{\mu_i} \Big|,
\end{align*}
where $\mu_i$ denotes the eigenvalue of $A$. One can check that
\[\sigma_{\min}( L_{\infty,h} ) \gtrsim 1, \qquad \min_i \Big|\frac{1-\e^{\mu_i t}}{\mu_i} \Big| \gtrsim \frac{1}{\|\bb{\lambda}\|_{\max} + \|\bb{\omega}\|/\omega_{\infty}},\]
where we have used the fact that $\|L_{\infty, h}\| \le 1/\omega_{\infty}$.
This completes the proof.
\end{proof}

The final solution can be obtained by solving the linear system $(I_h -\omega_{\infty}L_{\infty,h}) \bb{u}_h(T) = \bb{u}_0 + \bb{u}_T$.
According to Theorem \ref{thm:inputmodel}, we can apply the block-encoding oracle of $(I_h -\omega_{\infty}L_{\infty,h})^{-1}$ to \eqref{u0plusuT} and obtain
    \begin{align}
\ket{0} \otimes \ket{0^{n_a}} \otimes \ket{0^{n_s}}
&  ~~ ... ~
\xrightarrow{~R_t^\dag \otimes I^{\otimes n_a} \otimes I^{\otimes n_s}~}
\quad
\frac{1}{\eta_0 + \eta_1}\ket{0}\otimes \ket{0^{n_a}}
\otimes ( \bb{u}_0 + \bb{u}_T) + \ket{\bot}  \nonumber \\
& \qquad \xrightarrow{~~~~~~I \otimes \mathcal{U}_{\text{inv}} ~~~~~}
\quad
\frac{1}{(\eta_0 + \eta_1)\alpha_{\text{inv}}}\ket{0}\otimes \ket{0^{n_a}}
\otimes \bb{u}_h(T) + \ket{\bot}. \label{recover}
\end{align}

\subsection{Time complexity}

For the finite difference discretization, we can evaluate the error in the discrete $L^2$ norm, defined by
\[\|\bb{u}\|_h : =  \Big(h^d \sum_{\bb{i}} u_{\bb{i}}^2 \Big)^{1/2} = h^{d/2} \|\bb{u}\|, \]
where $\bb{i} = (i_1, \cdots,i_d)$ is the multi-index. In the following discussion, we assume homogeneous boundary conditions. For simplicity, we assess the complexity based on the queries to $\mathcal{U}_A$ in the following.

The goal of this section is to determine the computational cost of obtaining a quantum solution $\bb{u}_h^q$ that satisfies the following condition
\[\|\bb{u}(T) - \bb{u}_h^q(T)\|_h \le \delta = \delta(h).\]
The error can be divided into two parts:
\[ \|\bb{u}(T) - \bb{u}_h^q(T)\|_h  \le  \|\bb{u}(T) - \bb{u}_h(T)\|_h
+ \|\bb{u}_h(T) - \bb{u}_h^q(T)\|_h.\]
We will require each part of the error to be $\mathcal{O}(\delta)$, where $\bb{u}_h^q$ is obtained from \eqref{uh} with $\bb{U}_h(T)$ replaced by the approximation $\bb{U}_h^q(T)$.

The first term, which corresponds to the finite difference error, can be estimated using Lemma \ref{lem:errODE}:
\[\|\bb{u}(T) - \bb{u}_h(T)\|_h  \lesssim  \frac{\|\bb{\omega}\|}{\|\bb{\lambda}\|_{\min}} d h^2 =:\delta. \]
For the second term, we obtain from \eqref{utmodified} that
\begin{align*}
\|\bb{u}_h(T) - \bb{u}_h^q(T)\|_h
& \simeq h^{d/2} \| \bb{u}_h(T)- \bb{u}_h^q(T)\| \\
& = h^{d/2} \| ( I - \omega_{\infty} L_{h,d})^{-1} S (\bb{U}_h(T) - \bb{U}_h^q(T))  \| \\
& \lesssim h^{d/2}\|\bb{\omega}\|^{1/2} \|\bb{U}_h(T) - \bb{U}_h^q(T)\|,
\end{align*}
where we have used the fact that the matrix $L_h$ is negative definite.
Requiring the right-hand side to be $\mathcal{O}(\delta)$, we can take
\[\|\bb{U}_h(T) - \bb{U}_h^q(T)\| \lesssim  \frac{\delta}{h^{d/2}\|\bb{\omega}\|^{1/2}}
= \frac{\|\bb{\omega}\|^{1/2}}{\|\bb{\lambda}\|_{\min}} d h^{2-d/2} . \]

\begin{theorem}
Let $\delta = \frac{\|\bb{\omega}\|}{\|\bb{\lambda}\|_{\min}} d h^2$. Suppose the time-fractional heat equation \eqref{fracheat} is subject to homogeneous Dirichlet boundary conditions. If $T \lesssim \|\bb{\lambda}\|_{\max} + d h^{-2}(1 + \omega_{\infty} d h^{-2} )\|\bb{\omega}\|$, then there is a quantum algorithm that prepares an approximation of the normalized solution $\ket{\bb{u}_h(T)}$, denoted as $\ket{\bb{u}_h^q(T)}$, with $\Omega(1)$ probability and a flag indicating success, using $\widetilde{\mathcal{O}}(T^2 d^4 h^{-8})$ queries to the oracle $\mathcal{U}_A$. Here, $\widetilde{\mathcal{O}}$ suppresses all logarithmic factors, and multiplicative factors from $\bb{\lambda}$, $\bb{\omega}$ and $\bb{u}_0$ are omitted. The unnormalized vector $\bb{u}_h^q(T)$ provides a $\delta$-approximation of $\bb{u}_h(T)$ in the discrete $L^2$ norm.
\end{theorem}
\begin{proof}
We first consider the time complexity of solving the ODE system \eqref{systemmodifiedfinal}, with the unitary given by \eqref{Vunitary}.
According to Lemma \ref{lem:Schrcost}, we can prepare an $\varepsilon$-approximation, denoted by $\ket{\bb{U}^q_h(T)}$ of the normalized solution $\ket{\bb{U}_h(T)}$ with $\Omega(1)$ probability and a flag indicating success, using
\[\mathcal{C} = \widetilde{\mathcal{O}}\Big(\frac{T\|\bm{F}_h\|_{\text{smax}}}{\|\bb{U}_h(T)\|} \alpha_H  T \log^2 \frac{1}{\varepsilon}\Big)\]
queries to the oracle $\mathcal{U}_A$, where
\[\bm{F}_h = \bb{\omega}^{1/2} \otimes \bb{f}_h(t), \qquad
\bb{f}_h(t) = L_{h,d} \bb{u}_0.\]
By definition,
\[\|\bm{F}_h\|_{\text{smax}} \le \|\bb{\omega}\| \|L_{h,d} \bb{u}_0\| \lesssim  d h^{-2} \|\bb{\omega}\| \|\bb{u}_0\|.\]
The associate unitary operator is denoted by $V^q$, which can be written as
\[ \ket{0^{dn_x+m}} \quad \xrightarrow{ V^q } \quad  \ket{\bb{U}^q_h(T)}.\]

In the discussion of Section \ref{subsec:recover}, we should replace $V$ in \eqref{Vunitary} by $V^q $ as defined above.
Using the inequality $\| \frac{\bb{x}}{\|\bb{x}\|} - \frac{\bb{y}}{\|\bb{y}\|} \| \le 2 \frac{\|\bb{x} - \bb{y}\|}{\|\bb{x}\|}$ for two vectors $ \bb{x}, \bb{y} $, we can bound the error in the quantum state after a successful measurement as
\[
\| \ket{\bb{U}_h(T)} - \ket{\bb{U}_h^q(T)} \| \le \frac{2 \| \bb{U}_h(T) - \bb{U}^q_h(T) \|}{\| \bb{U}_h(T) \|}
\lesssim \frac{\|\bb{\omega}\|^{1/2}}{\|\bb{\lambda}\|_{\min}} \frac{d h^{2-d/2}}{\| \bb{U}_h(T) \|}=:\varepsilon.
\]
This along with Lemma \ref{lem:boundUh} indicates that we can take
\[
\frac{1}{\varepsilon} \lesssim \frac{\|\bb{\lambda}\|_{\min}}{\|\bb{\omega}\|^{1/2}} \frac{1}{d h^{2-d/2}} \frac{\|\bb{\omega}\|^{1/2}\|\bb{u}_0\|}{\|\bb{\lambda}\|_{\min}} = \frac{1}{d h^{2-d/2}} \|\bb{u}_0\| \lesssim h^{-2}.
\]

According the discussion in Section \ref{subsec:recover}, by applying amplitude amplification as described in \cite{BerryChilds2017ODE}, we can get $\Omega(1)$ probability with
\[\mathcal{O}((\eta_0 + \eta_1)\alpha_{\text{inv}} ) = \mathcal{O}( (\|\bb{u}_0\| + \|\bb{\omega}\|^{1/2} \|\bb{U}_h(T)\|) (1 + \omega_{\infty} d h^{-2} ) ) \]
repetitions of the above procedure. By applying Lemma \ref{lem:boundUh} once more and Eq.~\eqref{alphaA}, one can find the final complexity is
\begin{align*}
\mathcal{C}
& = \widetilde{\mathcal{O}}\Big(\frac{\|\bb{u}_0\| + \|\bb{\omega}\|^{1/2} \|\bb{U}_h(T)\|}{\|\bb{U}_h(T)\|} \|\bm{F}_h\|_{\text{smax}}\alpha_H  T^2 \log^2 \frac{1}{\varepsilon}(1 + \omega_{\infty} d h^{-2} )\Big) \\
& = \widetilde{\mathcal{O}}\Big(\frac{\|\bb{u}_0\| + \|\bb{\omega}\|^{1/2} \|\bb{U}_h(T)\|}{\|\bb{U}_h(T)\|} d h^{-2} \|\bb{\omega}\| \|\bb{u}_0\| \\
& \qquad \times (\|\bb{\lambda}\|_{\max} + \|\bb{\omega}\| dh^{-2}(1 + \omega_{\infty} d h^{-2}) +T)  T^2 \log^2 h^{-1} (1 + \omega_{\infty} d h^{-2} )\Big) \\
& = \widetilde{\mathcal{O}}\Big(\frac{\|\bb{u}_0\| + \|\bb{\omega}\|^{1/2} \|\bb{U}_h(T)\|}{\|\bb{U}_h(T)\|} d^2 h^{-4} \|\bb{\omega}\|^2 \|\bb{u}_0\|\|\bb{\lambda}\|_{\max}  T^2 (1 + \omega_{\infty} d h^{-2} )^2\Big) \\
& = \widetilde{\mathcal{O}}\Big(( \frac{\|\bb{\lambda}\|_{\max} + \|\bb{\omega}\|/\omega_{\infty}}{\|\bb{\omega}\|^{1/2}} + \|\bb{\omega}\|^{1/2}) d^2 h^{-4} \|\bb{\omega}\|^2 \|\bb{u}_0\|\|\bb{\lambda}\|_{\max}  T^2 (1 + \omega_{\infty} d h^{-2} )^2 \Big) \\
& = \widetilde{\mathcal{O}}\Big(T^2 d^2 h^{-4}(1 + \omega_{\infty} d h^{-2} )^2 \|\bb{\omega}\|^{3/2} \|\bb{\lambda}\|_{\max} (\|\bb{\lambda}\|_{\max} + \|\bb{\omega}\|/\omega_{\infty}) \|\bb{u}_0\|\Big),
\end{align*}
where we have used the assumption $T \lesssim \|\bb{\lambda}\|_{\max} + d h^{-2}(1 + \omega_{\infty} d h^{-2} )\|\bb{\omega}\|$.
\end{proof}

 To solve \eqref{systemmodifiedfinal} or \eqref{ukeq} using the classical method, we consider the forward Euler discretization, given by
\begin{align}
&  \frac{\tilde{\bb{u}}_{h,n+1}^{(k)} - \tilde{\bb{u}}_{h,n}^{(k)}}{\Delta t} + \lambda_k \tilde{\bb{u}}_{h,n}^{(k)}= (1+\lambda_k)   \sum_{i=1}^M \frac{\omega_i}{1+\lambda_i}  L_{\infty,h}  \tilde{\bb{u}}_{h,n}^{(i)}   + (1+\lambda_k) \bb{f}_{h,n},   \label{ukeqtime} \\
& \tilde{\bb{u}}_{h,n}^{(k)}(0) = \bb{0}, \qquad k = 1,\cdots, M, \label{ukeq0time}
\end{align}
where $\tilde{\bb{u}}_{h,n}^{(i)}$ is the approximation of $\tilde{\bb{u}}_{h}^{(i)}(t_n)$.
Noting that $ L_{\infty,h} = L_{h,d} ( I_h - \omega_{\infty} L_{h,d})^{-1} $, at each time step, we need to compute
\[
L_{\infty,h} \tilde{\bb{u}}_{h,n}^{(i)} = L_{h,d} ( I_h - \omega_{\infty} L_{h,d})^{-1} \tilde{\bb{u}}_{h,n}^{(i)}, \quad i = 1, \cdots, M.
\]
As shown in the proof of Theorem 3.1 in \cite{JLY2022multiscale}, the conjugate gradient method can be applied to solve the linear system $ ( I_h - \omega_{\infty} L_{h,d}) \bb{x}^{(i)} = \tilde{\bb{u}}_{h,n}^{(i)} $, with the time complexity given by
\[
\mathcal{O}(N s \sqrt{\kappa}) = \mathcal{O}(h^{-d} d \sqrt{h^{-1}}) = \mathcal{O}(d h^{-(d+0.5)}),
\]
where $ N $, $ s $, and $ \kappa $ represent the number of entries in each row or column, the sparsity, and the condition number of the underlying matrix, respectively. The number of basic operations involved on computing $L_{h,d} \bb{x}^{(i)}$ is $\mathcal{O}(N s) = \mathcal{O}(d h^{-d})$. Therefore, at each time step, we use at least $\mathcal{O}(M d h^{-(d+0.5)})$ matrix-vector multiplications for the forward Euler approach.

According to Lemma \ref{lem:errODE}, the error in the discrete $ L^2 $ norm is $ \mathcal{O}(\tau + d h^{-2}) $. Given that the stability condition for the time discretization is $ d \tau/h^2 \le 1 $, where $ \tau = T / N_t $, we can choose $ \tau = \mathcal{O}(h^2 / d) $ or $ N_t = \mathcal{O}(T d h^{-2}) $.
In contrast, the quantum algorithm exhibits a dependence on $h^{-1}$ that is {\it independent} of the dimensionality, while the classical method scales exponentially, particularly in high dimensions. This highlights the exponential quantum advantage for high-dimensional fractional problems, where classical methods are hindered by the curse of dimensionality.

\section{Numerical simulation} \label{sec:numerical}

We consider the $d$-dimensional fractional heat equation in \cite{Khristenko23fractional} with homogeneous Dirichlet boundary conditions:
\[\begin{cases}
\partial_t^{\alpha} u(t,x) - \Delta u(t,x) = 0, \qquad & t\in (0,T), ~~x = (x_1,\cdots,x_d)\in \Omega = (0,1)^d, \\
u(0,x) = \sin(\pi x),  \qquad & x \in [0,1]^d,\\
u(t,x)|_{\partial \Omega} = 0, \qquad & t\in [0,T],
\end{cases}\]
where $\sin(\pi x) := \sin(\pi x_1)\cdots \sin(\pi x_d)$.
The analytical solution is given by
\[u(t,x) = E_{\alpha,1}[-d \pi^2 t^{\alpha}] \sin (\pi x), \]
where
\[E_{\alpha, \beta}[z] = \sum_{k=0}^{(\infty)} \frac{z^k}{\Gamma(\alpha k + \beta)}\]
is the Mittag-Leffler function.

\subsection{The finite difference discretization in 1D}

\begin{figure}[!htb]
		\centering
        \includegraphics[scale=0.5]{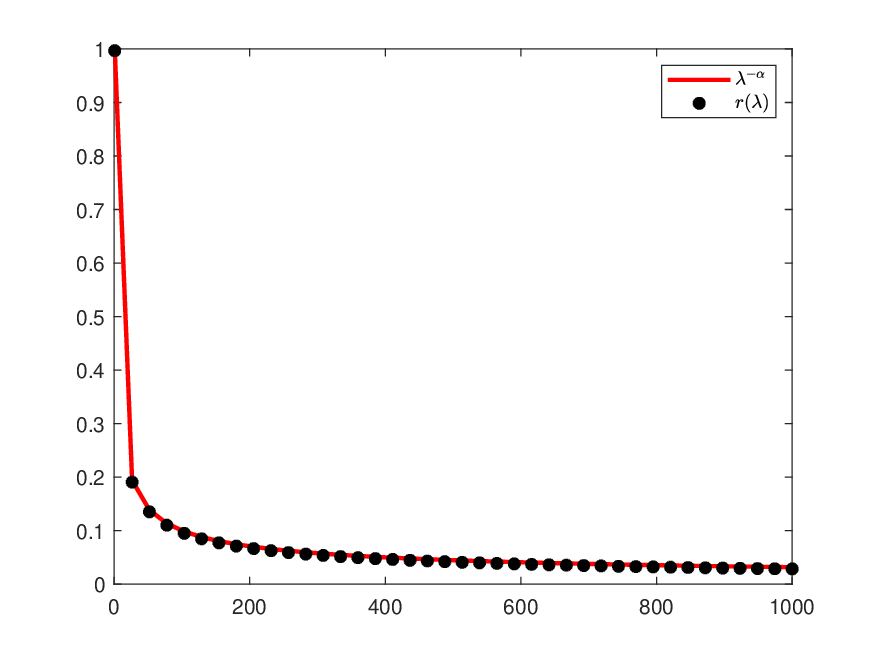}
		\caption{Rational approximation of $\lambda^{-\alpha}$}\label{fig:ratapp}
	\end{figure}

For the rational approximation of $\lambda^{-\alpha}$, we apply the AAA algorithm \cite{Nakatsukasa2018AAA}. The implementation is carried out using the built-in function `aaa.m' from the open-source Matlab package, {\it Chebfun}. We set the AAA tolerance to $10^{-6}$ with 1000 candidate points.
The parameter in \eqref{lambdafrac} is taken as $\tau = 0.001$.
The approximation for $\alpha = 0.5$ is shown in Fig.~\ref{fig:ratapp}.

For the spatial discretization, we choose $N_x-1 = 2^5$.
For the setting of the Schr\"odingerization approach, please refer to \cite{JLMY2025SchOptimal}.

\begin{figure}[!htb]
		\centering
        \subfigure[$T=1$]{\includegraphics[scale=0.45]{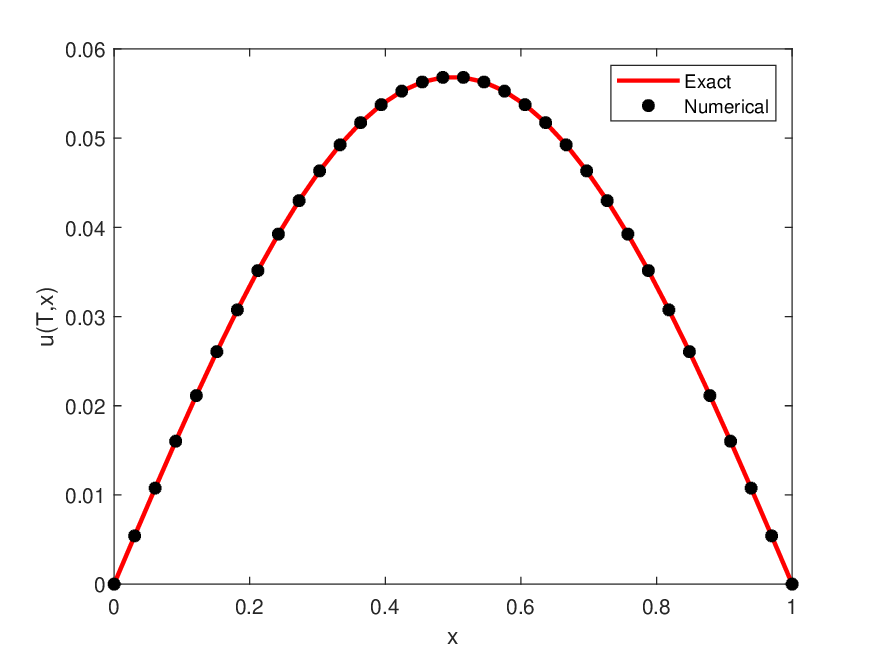}}
		\subfigure[$T=2$]{\includegraphics[scale=0.45]{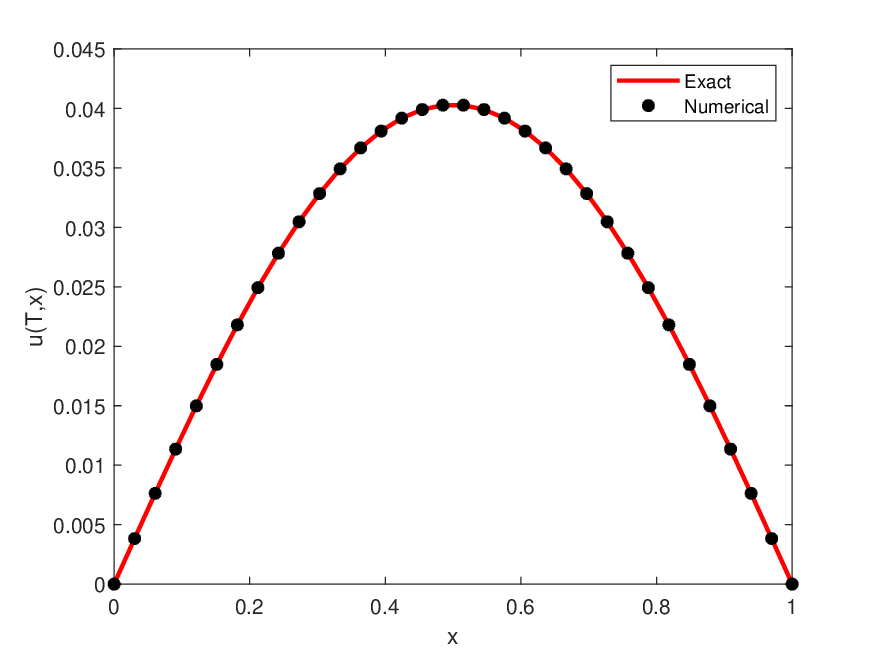}}\\
		\caption{Numerical and exact solutions for $\alpha = 0.5$}\label{fig:Tfrac}
	\end{figure}

The numerical and exact solutions for $\alpha = 0.5$ are displayed in Fig.~\ref{fig:Tfrac},
from which we observe that the Schr\"odingerization based approach produces accurate approximate results.

\subsection{The finite element discretization in 2D}

Given a triangular mesh $\mathcal{T}_h$ for the spatial domain $\Omega = (0,1)^2$, we define the continuous piecewise linear finite element space on the mesh as
\begin{equation}\label{Vl}
	V_h:=\{v\in H_0^1(\Omega): v|_K \in \mathbb{P}_1(K), \quad K\in \mathcal{T}_h \},
\end{equation}
where $\mathbb{P}_1$ is the space of linear polynomials. Let $V_h = \text{span}\{ \varphi_i: i = 1,\cdots, N_h\}$ with $\varphi_i$ being the nodal basis function and $N_h$ the number of the basis functions. In the implementation, we consider the tensor-product linear finite element.

	\begin{figure}[!htb]
		\centering
        \subfigure[Exact solution]{\includegraphics[scale=0.45]{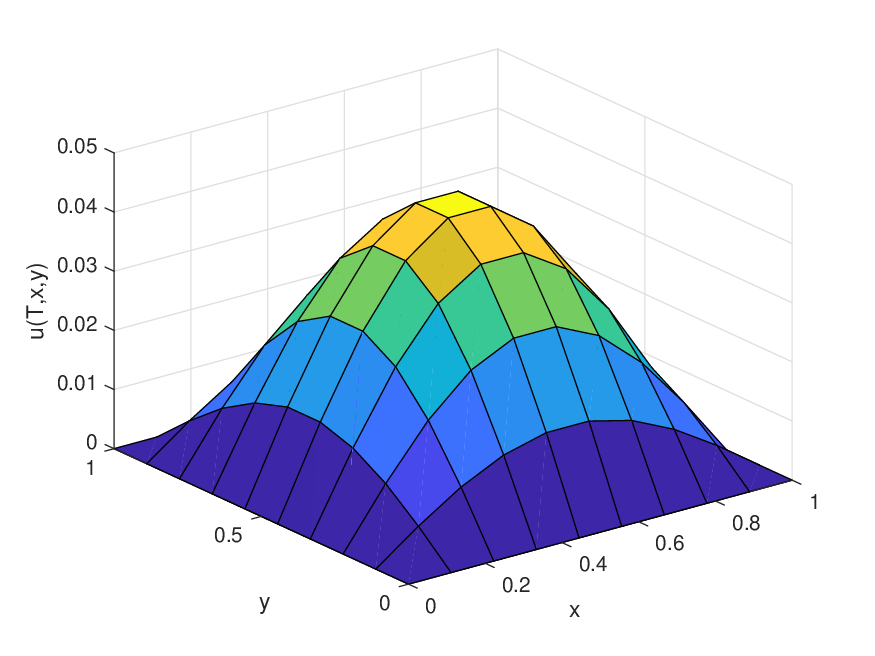}}
		\subfigure[FEM solution]{\includegraphics[scale=0.45]{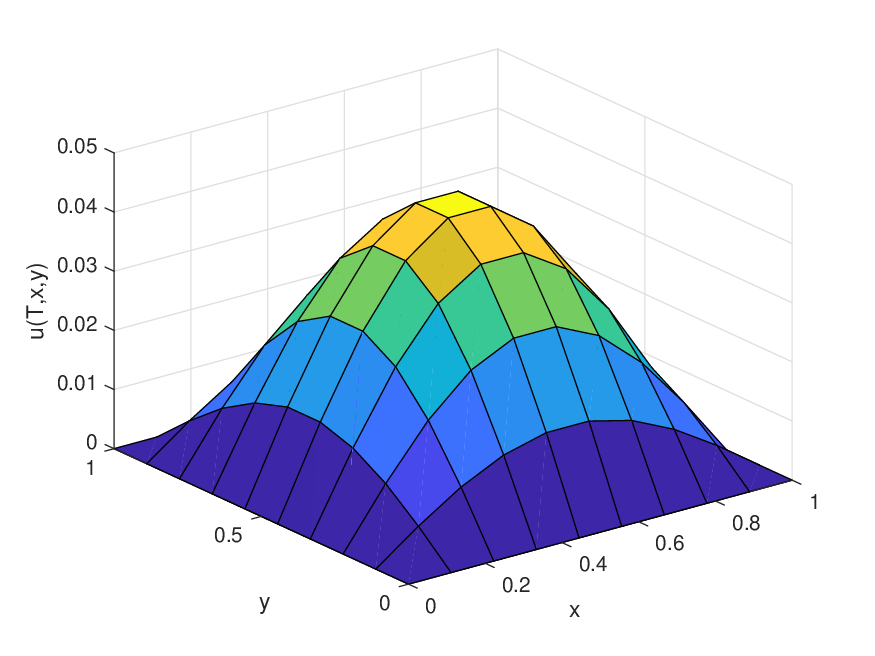}}\\
		\caption{Numerical and exact solutions for $\alpha = 0.1$ with FEM discretization ($T=1$)}\label{fig:TfracFEM}
	\end{figure}

The spatial discretization of \eqref{uk}-\eqref{uinf} is then given by
\begin{align*}
& \frac{1}{1+\lambda_k} \frac{\d }{\d t} M_h \tilde{\bb{u}}^{(k)}(t) + \frac{\lambda_k}{1+\lambda_k} M_h \tilde{\bb{u}}^{(k)}(t) - A_h \bb{u} = \bb{0}, \qquad k = 1,\cdots, M,  \\
& M_h\tilde{\bb{u}}^{(\infty)}(t) -A_h \bb{u}(t) = \bb{0},
\end{align*}
where $M_h = ((\varphi_i, \varphi_j))$ and $A_h = ((\nabla \varphi_i, \nabla \varphi_j))$ are the mass and stiffness matrices, respectively. In this case, one can express \eqref{uexpand} as
\[
\bb{u}(t) =  \bb{u}_0 + \sum_{i=1}^M \frac{\omega_i}{1+\lambda_i} \tilde{\bb{u}}^{(i)}(t) + \omega_{\infty} \tilde{\bb{u}}^{(\infty)}(t).
\]
We can solve the system by applying a similar procedure for the finite difference discretization. In fact, one can treat $ M_h \tilde{\bb{u}} $ as the vector $ \tilde{\bb{u}} $ in \eqref{schemeFDM}.
The finite element solution is obtained via
\[u_h(t) = \sum_{i=1}^{N_h} \bb{u}_i \varphi_i(x),\]
where $\bb{u}_i$ is the $i$-th entry of $\bb{u}$.

We plot the numerical and exact solutions for $\alpha = 0.1$ in Fig.~\ref{fig:TfracFEM}, where we observe that the Schr\"odingerization based approach yields accurate approximate results, thereby validating the correctness of the dimension-lifting approach.

\section{ Conclusion}

In this article we developed a quantum algorithm for simulating high-dimensional time-fractional heat equations.
The algorithm combines two key techniques: (1) the dimension-lifting technique introduced in \cite{Eritz23timeFractional}, which converts the nonlocal fractional gradient flow to an equivalent integer-order gradient flow in one higher dimension, and (2) the Schr\"odingerization method, which converts the resulting PDE with finite difference or finite element discretizations into a Hamiltonian simulation problem.

We performed a comprehensive analysis of the time complexity, demonstrating that the quantum algorithm offers an exponential advantage over classical methods in high-dimensional settings. The quantum algorithm achieves a query complexity of $ \widetilde{\mathcal{O}}(T^2 d^4 h^{-8}) $, where the scaling in nverse mesh size $ h^{-1} $ is independent of the dimension $ d $. In contrast, classical approaches require at least $ \mathcal{O}(d h^{-(d+0.5)}) $ operations per time step, emphasizing the potential of quantum computing to mitigate the curse of dimensionality in fractional PDEs.

Numerical experiments validated the correctness of our method, highlighting the versatility of the Schr\"odingerization framework. This framework can be adapted to other linear time-fractional problems, such as the time-fractional biharmonic equation, as well as equations with physical boundary conditions. Additionally, quantum circuit implementations based on \cite{Sato24Circuit, JLY24Circuits, HuJin24SchrCircuit} warrant further exploration.

\section*{Acknowledgments}

Shi Jin thanks Prof. Barbara Wahlmuth for pointing out reference \cite{Eritz23timeFractional}, which provides the starting point of this article.

SJ and NL are supported by NSFC grant No. 12341104,
the Shanghai Jiao Tong University 2030 Initiative and the Fundamental Research Funds for the Central Universities. SJ was also partially supported by the NSFC grant No. 12426637 and  the Shanghai Municipal Science and Technology Major Project (2021SHZDZX0102). NL also acknowledges funding from the Science and Technology Commission of Shanghai Municipality (STCSM) grant no. 24LZ1401200 (21JC1402900), NSFC grant No.12471411 and the Shanghai Science and Technology Innovation Action Plan (24LZ1401200).
YY was supported by NSFC grant No. 12301561, the Key Project of Scientific Research Project of Hunan Provincial Department of Education (No. 24A0100) and the 111 Project (No. D23017).

\bibliographystyle{alpha} 
\bibliography{Refs}

\end{document}